\documentclass[onefignum,onetabnum]{siamart190516}

\usepackage{amsfonts}
\usepackage{amsmath}
\usepackage{amssymb}
\usepackage{bm}
\usepackage{graphicx}
\usepackage{epstopdf}
\usepackage{caption}
\usepackage{subcaption}
\usepackage{algorithm}
\usepackage{algpseudocode}
\usepackage{setspace}
\usepackage{ulem}
\usepackage{multirow}
\usepackage{booktabs}
\ifpdf
  \DeclareGraphicsExtensions{.eps,.pdf,.png,.jpg}
\else
  \DeclareGraphicsExtensions{.eps}
\fi
   
\algnewcommand{\IIf}[1]{\State\algorithmicif\ #1\ \algorithmicthen}

\usepackage[version=4]{mhchem}

\usepackage{epsfig, algorithm, algpseudocode}
\usepackage{multirow, booktabs}
\usepackage{xr}

\newsiamremark{remark}{Remark}
\newsiamremark{hypothesis}{Hypothesis}
\crefname{hypothesis}{Hypothesis}{Hypotheses}
\newsiamthm{claim}{Claim}

\title{Partition of Unity Data-driven Finite Element Method for Anisotropic Diffusion in Magnetized Plasmas}

\author{
  Siu Wun Cheung\thanks{Center for Applied Scientific Computing, Lawrence Livermore National Laboratory, Livermore, CA 94550 (\email{cheung26@llnl.gov}, \email{choi15@llnl.gov}} \and
  Youngsoo Choi\footnotemark[1]
}

\begin{document}

\maketitle

\begin{abstract}
Numerical simulation of transport processes in magnetized plasmas is characterized by extreme anisotropy, where diffusion along magnetic field lines can be many orders of magnitude larger than in the perpendicular direction. This anisotropy leads to thin, steep boundary layers and numerical pollution that standard finite element methods struggle to resolve without prohibitive mesh refinement.
In this paper, we introduce the Partition of Unity Data-Driven Finite Element Method (PU-DDFEM), a data-driven local model reduction approach, for anisotropic diffusion. 
The method utilizes a \textit{bottom-up} construction to enable a complete offline-online splitting. 
We parametrize the local magnetic field and boundary conditions to generate a snapshot space that captures the essential physics of plasma confinement. 
Dominant modes are extracted via a local spectral problem to identify representative features. 
These modes are then mapped to the global domain using a partition of unity to form a conforming approximation space. 
This ensures that the basis functions are fully reusable across different global magnetic configurations, significantly reducing the online computational overhead for large-scale plasma simulations.
\end{abstract}

\section{Introduction}
\label{sec:intro}

Achieving sustainable fusion requires a precise understanding of transport processes within the tokamak, where charged particles are confined by intense magnetic fields. A defining characteristic of these systems is the extreme anisotropy of the diffusion tensor. As particles follow the helical trajectories of magnetic field lines, transport in the direction parallel to the field is remarkably efficient, while transport in the perpendicular direction is suppressed by several orders of magnitude. The physics of magnetically confined fusion plasmas is characterized by an extreme anisotropy in transport with respect to the direction of the magnetic field. Due to their strong magnetization, charged particles follow helical trajectories around magnetic field lines, resulting in significantly faster transport along the magnetic field than across it \cite{braginskii1965transport}. This anisotropy is particularly pronounced in thermal conduction, where the parallel Spitzer--Härm conductivity, dominated by electron transport, can exceed its perpendicular counterpart by several orders of magnitude \cite{spitzer1953transport}. In reactor-relevant tokamak regimes, the ratio of parallel to perpendicular thermal conductivity, $\kappa_\parallel/\kappa_\perp$, can reach values as high as $10^9$--$10^{12}$, particularly near the plasma edge \cite{gunter2005modelling,vogl2023mesh}. Such extreme anisotropy has important consequences for the spatial structure of plasma transport. In particular, the interaction between the anisotropic transport and the boundary conditions at the plasma edge gives rise to narrow boundary layers, commonly associated with the scrape-off layer. The characteristic width of these layers scales approximately as $\sqrt{\kappa_\perp/\kappa_\parallel}$ relative to the characteristic domain size, making them exceedingly thin even when perpendicular transport is enhanced by turbulent processes.

The numerical simulation of magnetically confined fusion plasmas remains one of the most demanding challenges in computational physics. From a computational perspective, resolving the fine-scale structures presents a major challenge for conventional numerical discretizations. In particular, standard finite element methods can suffer from severe numerical pollution when the computational mesh is not sufficiently aligned with the magnetic field. Discretization errors may introduce artificial perpendicular transport through the leakage of the much larger parallel flux into the perpendicular direction. Consequently, accurate simulations often require specialized numerical treatments to capture the local physics, such as highly anisotropic meshes aligned with the magnetic field, carefully designed discretization schemes for non-aligned coordinates \cite{gunter2005modelling} or extensive local mesh refinement \cite{vogl2023mesh}.

In this paper, we propose a data-driven approach to tackle this computational challenge. We follow the principle of data-driven finite element method (DD-FEM) introduced in \cite{choi2025defining} and propose a bottom-up pre-trained local model reduction approach designed to circumvent these resolution requirements by capturing the local physics of the plasma within a compressed approximation space. Localized model reduction is uniquely suited for anisotropic diffusion because the most difficult numerical features, oscillatory gradients and boundary layers, are fundamentally driven by the local orientation and magnitude of the magnetic field. By pre-solving specialized local problems that learn these features, we can construct a coarse-scale discretization that remains accurate even when the coarse mesh does not explicitly resolve the fine-scale layers.

Within the DD-FEM architecture, these data-driven elements are built by leveraging component reduced order models \cite{huynh2013static, huynh2013static2, eftang2014port, mcbane2021component, mcbane2022stress, chung2024train, chung2024scaled, chung2024scalable, zanardi2024scalable, hoang2021domain, ebrahimi2024hyperreduced, smetana2016optimal, eftang2013port, diaz2024fast, wentland2025role, vasilyeva2021preconditioning}. This strategy relies on compressing fine-scale simulation records over discretized element domains to form compact bases, subsequently projecting the local governing equations onto these low-dimensional subspaces. Grounded in standard projection-based model reduction principles \cite{benner2015survey, choi2019space, carlberg2018conservative, hesthaven2016certified, kim2022fast, lee2020model, quarteroni2015reduced, wang2012proper, tsai2025local}, this methodology guarantees both computational speedups and structural fidelity at the element level, supported by rigorous error bounds established in prior literature \cite{huynh2013static, eftang2014port, mcbane2021component}.

The global assembly approach used in this work is related to the broader line of work on Partition-of-Unity Finite Element Methods (PUFEM) and Generalized Finite Element Methods (GFEM) \cite{melenk1996partition, babuska1997partition, babuvska2004generalized}, including multiscale extensions such as Generalized Multiscale Finite Element Methods (GMsFEM) \cite{efendiev2013generalized, calo2016randomized}. These methods provide a natural framework for enriching a coarse finite element space with localized functions that encode unresolved features of the solution. The main difference of our approach to previous works lies in how the local approximation spaces are constructed and deployed via a bottom-up training strategy. An important contribution of this work is the development of a systematic geometric parametrization that captures the spatially varying orientation, local curvature, and extreme anisotropy ratios of the magnetic field. Rather than relying on problem-specific enrichment, the local magnetic field topologies and boundary data are treated as parameters within this geometric framework. A diverse collection of local problems is solved offline to generate a reusable library of snapshots, ensuring that the resulting approximation space can represent a broad family of local transport behaviors across varied magnetic configurations.

The bottom-up pre-training and global assembly strategy ensures a complete offline-online splitting. The basis functions are designed to be reusable, so once the offline training and parameterization are complete, the resulting subspace can be deployed to efficiently solve forward simulation problems across multiple static magnetic configurations without recomputing the local spectral modes. This makes the method particularly attractive for parametric studies and scenario exploration where the magnetic equilibrium may vary. This hybridization of data-driven local model reduction and variational approximation extends beyond standard anisotropic diffusion and provides a promising direction for interpretable, generalizable, and physically grounded data-driven discretizations.
\section{Steady-State Plasma Transport}
\label{sec:problem}

We consider the steady-state anisotropic diffusion equation as a model for heat and particle transport in magnetically confined plasmas, such as those encountered in tokamak devices. 
Following the mathematical formulation in \cite{vogl2023mesh}, we consider a bounded Lipschitz domain $\Omega\subset\mathbb{R}^2$ and seek a scalar field $u:\Omega\to\mathbb{R}$ satisfying
\begin{equation}
-\nabla\cdot\left(\mathbf{K}\nabla u\right)=f
\quad\text{in }\Omega,
\label{eq:strong_form}
\end{equation}
where $f\in L^2(\Omega)$ represents a volumetric heat or particle source, and $\mathbf{K} \in L^\infty(\Omega; \mathbb{S}_2^{++})$ is the spatially varying diffusion tensor induced by the magnetic field. Here, $\mathbb{S}_2^{++}$ denotes the open convex cone of all $2 \times 2$ symmetric and positive definite matrices. We also assume that $\mathbf{K}$ is uniformly positive definite. More precisely, there exist constants $0<k_{\min}\leq k_{\max}<\infty$ such that
\begin{equation}
k_{\min}|\boldsymbol{\xi}|^2
\leq
\boldsymbol{\xi}^{\top}\mathbf{K}(\mathbf{x})\boldsymbol{\xi}
\leq
k_{\max}|\boldsymbol{\xi}|^2
\quad
\forall\boldsymbol{\xi}\in\mathbb{R}^2
\end{equation}
for almost every $\mathbf{x} = (x,y)\in\Omega$. The ratio $k_{\max}/k_{\min}$ may be extremely large, reflecting the strong anisotropy of magnetized transport.
In this work, we consider benchmarks where the diffusion equation \eqref{eq:strong_form} is subject to 
homogeneous natural boundary conditions and possibly inhomogeneous essential boundary conditions that include prescribed Dirichlet data and periodic constraints. 
The specific boundary configuration depends on the magnetic-field geometry and will be described in the numerical experiments.
To accommodate these cases in a unified formulation, let $\mathcal{C}: H^1(\Omega) \to W(\partial \Omega)$ 
denote the operator specifying the essential boundary constraints, and let $g \in W(\partial \Omega)$ denote the associated boundary data. 
To formulate this, we introduce a lifting function $u_g \in H^1(\Omega)$ that satisfies the boundary constraint, such that $\mathcal{C}u_g = g$.
We define the test space by 
\begin{equation}
V_0
=\{ v \in H^1(\Omega): \mathcal{C}v = 0\}.
\end{equation}
The trial space is then defined as the affine space
\begin{equation}
V_g = \{ v \in H^1(\Omega) : \mathcal{C}v = g \} = u_g + V_0,
\end{equation}
where $u_g \in H^1(\Omega)$ is a lifting function  that satisfies the boundary constraint $\mathcal{C}u_g = g$. 
The weak formulation is to find $u \in V_g$ such that
\begin{equation}
a(u,v)=\ell(v)
\qquad\forall v\in V_0,
\end{equation}
or, equivalently, fund $u_0 \in V_0$ such that 
\begin{equation}
a(u_0,v) = \ell(v) - a(u_g, v)
\qquad\forall v\in V_0,
\label{eq:weak_form}
\end{equation}
where the bilinear form $a: H^1(\Omega) \times H^1(\Omega) \to \mathbb{R}$ and the linear form $\ell: H^1(\Omega) \to \mathbb{R}$ are defined by
\begin{equation}
\begin{split}
a(u,v) & = \int_\Omega
(\nabla v)^\top \mathbf{K} \nabla u \, \mathrm{d}\mathbf{x} \quad \forall u,v \in H^1(\Omega), \\
\ell(v) & =\int_\Omega vf \, \mathrm{d}\mathbf{x} \quad \forall v \in H^1(\Omega). \\
\end{split}
\end{equation}

\subsection{Magnetic Field Configuration}
\label{subsec:magnetic_field}

In a two-dimensional poloidal cross-section of a fusion plasma, the magnetic field geometry can be described using a scalar magnetic flux function $A_z$. The poloidal magnetic field is obtained from the curl of the vector potential: 
\begin{equation}
\mathbf{B}_p(\mathbf{x})
=
\nabla\times(A_z(\mathbf{x}) \mathbf{e}_3)
=
\begin{pmatrix}
\partial_y A_z(\mathbf{x})\\
-\partial_x A_z(\mathbf{x})\\
0
\end{pmatrix}.
\label{eq:poloidal_field}
\end{equation}
Consequently, the contours of $A_z$ describe the poloidal magnetic flux surfaces, while critical points of $A_z$ correspond to locations where the poloidal field vanishes.
In a tokamak, a nonzero toroidal component can regularize the direction of the full three-dimensional magnetic field. We therefore introduce a constant toroidal field
\begin{equation}
\mathbf{B}_t=B_0\mathbf{e}_3,
\qquad B_0>0,
\end{equation}
and define the total field by $\mathbf{B}=\mathbf{B}_p+\mathbf{B}_t$, where 
Its projection onto the poloidal plane is
\begin{equation}
\mathbf{b}(\mathbf{x})
=
\frac{\mathbf{B}_p(\mathbf{x})}
{\sqrt{|\mathbf{B}_p(\mathbf{x})|^2+B_0^2}},
\label{eq:projected_field}
\end{equation}
satisfying $0 \leq | \mathbf{b} |<1$ and vanishing at a poloidal magnetic null.

We construct the effective two-dimensional diffusion tensor by restricting the three-dimensional anisotropic transport tensor to the poloidal plane. Specifically, let $\kappa_{\parallel}>0$ and $\kappa_{\perp}>0$ denote the parallel and perpendicular diffusion coefficients, respectively. The resulting tensor is
\begin{equation}
\mathbf{K}(\mathbf{x})
=
\kappa_{\perp}\mathbf{I}
+
(\kappa_{\parallel}-\kappa_{\perp})
\mathbf{b}(\mathbf{x})\otimes\mathbf{b}(\mathbf{x}),
\label{eq:diffusion_tensor}
\end{equation}
where $\mathbf{I}$ is the $2\times2$ identity matrix and $\otimes$ denotes the outer product. This expression accounts for the magnitude of the poloidal projection of the full magnetic-field direction. In particular, the effective parallel enhancement decreases as the poloidal component becomes small, while the perpendicular contribution remains positive.
The anisotropy ratio quantifies the disparity between the transport coefficients. When $\kappa_{\parallel} \gg \kappa_{\perp}$, transport along the magnetic field is much more efficient than cross-field transport. In the limiting regime of very large anisotropy, the solution tends toward a field-aligned equilibrium in regions where the field geometry and boundary conditions permit such behavior. Sharp gradients can develop across magnetic separatrices and other interfaces separating distinct field-line regions.
This structure motivates a local model reduction strategy in which the dependence of the solution on the local field orientation, spatial variation, and anisotropy is represented through precomputed approximation spaces. Rather than attempting to resolve every thin layer using global mesh refinement alone, the proposed method enriches a coarse-scale approximation with local modes trained to represent characteristic anisotropic transport patterns.

\subsection{Finite Element Solution}
\label{subsec:fem_solution}

We first introduce a high-fidelity finite element discretization that serves as the reference problem for constructing and evaluating the reduced approximation. Let $\mathcal{T}_h$ be a conforming mesh of $\Omega$, with characteristic mesh size $h$, and let $V_h \subset H^1(\Omega)$ denote a $H^1$-conforming finite element space, and $V_{h,0} = V_h \cap V_0$ be its subspace subject to the homogeneous essential boundary conditions. The finite element solution $u_h = u_{h,g} + u_{h,0} \in V_h$ of \eqref{eq:strong_form} is obtained by taking $u_{h,g} = \mathcal{J}_{\Omega} u_g \in V_h$ where $\mathcal{J}_{\Omega}$ is the global interpolation operator onto $V_{h}$, and solving $u_{h,0} \in V_{h,0}$ such that 
\begin{equation}
a(u_{h,0},v)
= \ell(v) - a(u_{h,g},v)
\quad\forall v\in V_{h,0}.
\label{eq:fem}
\end{equation}
Let $\{v_j\}_{j=1}^{N_h}$ denote the standard shape basis functions of $V_{h,0}$, where $N_h = \dim(V_{h,0})$. Expressing the finite element solution as $u_{h,0} = \sum_{j=1}^{N_h} U_j v_j$ and choosing test functions $v = v_i$ for $i = 1, \dots, N_h$ transforms \eqref{eq:fem} into the sparse linear algebraic system on $\mathbf{U}_h = (U_j)_{j=1}^{N_h}  \in \mathbb{R}^{N_h}$ 
\begin{equation}
\mathbf{A}_h \mathbf{U}_h = \mathbf{F}_h,
\label{eq:full_linear}
\end{equation}
where the global stiffness matrix $\mathbf{A}_h = (A_{ij})_{i,j=1}^{N_h} \in \mathbb{R}^{N_h \times N_h}$ and the load vector $\mathbf{F}_h = (F_i)_{i=1}^{N_h} \in \mathbb{R}^{N_h}$ are assembled by summing the local element contributions, with entries defined by
\begin{equation}
A_{ij} = a(v_j, v_i), \quad F_i = \ell(v_i) - a(u_{h,g},v_i), \quad i,j = 1, \dots, N_h.
\end{equation}
Thanks to the sparsity, symmetry, and positive definiteness of the matrix $\mathbf{A}$, the linear system \eqref{eq:full_linear} is solved using the preconditioned conjugate gradient (PCG) method equipped with an algebraic multigrid (AMG) preconditioner.
Despite the well-posedness inheritied from the weak form \eqref{eq:weak_form}, the numerical approximation can become challenging when $\kappa_{\parallel} / \kappa_{\perp}$ is large. As discussed by Vogl et al.~\cite{vogl2023mesh}, a mesh that is not sufficiently aligned with the magnetic field can produce artificial cross-field transport. The large parallel diffusion coefficient amplifies the effect of discretization errors, so that the resulting numerical pollution may dominate the physical perpendicular diffusion. Resolving the resulting layers using conventional finite elements can therefore require substantial mesh refinement.
The difficulty is not simply the number of degrees of freedom: it is also the need to represent solution features whose orientation and spatial scale are determined by the local magnetic geometry. The PU-DDFEM framework introduced in the next section addresses this difficulty by constructing localized approximation functions that encode representative anisotropic transport behavior. The resulting reduced Galerkin problem retains the variational structure of the weak form \eqref{eq:weak_form}, while seeking an accurate approximation in a substantially smaller space.

\section{Partition-of-Unity Data-Driven Finite Element Method}
\label{sec:lmr}

The objective of the family of data-driven finite element methods is to construct a low-dimensional approximation space that captures the dominant local features of the underlying physics while reducing the computational cost of repeated global solves. The central idea is to separate the construction of local approximation spaces from their deployment on a particular global problem. 
In this work, we consider the partition-of-unity formulation, where the domain $\Omega$  
is divided into overlapping patches $\left\{\omega_i\right\}_{i \in \mathcal{I}}$.
We assume that the mesh $\mathcal{T}_h$ is conforming to the cover $\left\{\omega_i\right\}_{i \in \mathcal{I}}$, and 
there is a reference patch $\widehat{\omega}$ which can be  
mapped to each physical patch $\omega_i$ via a smooth bijection  
$T_i:\omega_i\to\widehat{\omega}$.

In the context of anisotropic diffusion in magnetized plasmas, we parameterize local magnetic configurations, generate fine-scale local snapshots, and compress the resulting snapshot spaces into a library of reduced bases in. This bottom-up construction is intended to make the local bases reusable across a family of global magnetic configurations, provided that their local tensor fields are adequately represented by the offline parameterization. In the online phase, we assemble the selected local modes into a conforming approximation space $V_{H,0}$ using a partition of unity for solving the global problem. 

\subsection{Bottom-Up Pre-Training}
\label{subsec:offline}

The offline phase builds a library of local reduced-order models on a fixed reference patch. The training parameters describe representative variations in local magnetic-field orientation and magnitude, while the boundary excitations probe the response of the local diffusion operator to different boundary profiles. A subsequent local spectral reduction compresses these responses into a small set of basis functions associated with each region of the parameter space.

\subsubsection{Parametrization}
\label{subsubsec:parametrization}

The local magnetic configuration is parametrized by a five-dimensional parameter vector
\begin{equation}
\boldsymbol{\mu}_{\mathbf{B}}
=
(\varphi_0,\alpha_1,\alpha_2,\beta_1, \beta_2)
\in\mathcal{P}_{\mathbf{B}}\subset\mathbb{R}^5.
\label{eq:local_parameter}
\end{equation}
Here, $\varphi_0$ specifies its nominal orientation, 
$\alpha_1$ and $\alpha_2$ control the spatial variation of that orientation, and 
$\beta_1$ and $\beta_2$ control spatial variations in the magnitude of the poloidal projection.
We represent the magnitude profile $\widehat{b}$ and orientation angle $\widehat{\varphi}$ by: 
for all $\widehat{\mathbf{x}} = (\widehat{x}, \widehat{y}) \in \widehat{\omega}$,
\begin{equation}
\begin{aligned}
\widehat{\varphi}(\widehat{\mathbf{x}};\boldsymbol{\mu}_{\mathbf{B}})
&=
\varphi_0+\alpha_1 \widehat{x}+\alpha_2 \widehat{y}, \\
\widehat{b}(\widehat{\mathbf{x}};\boldsymbol{\mu}_{\mathbf{B}})
&=
\left(1+\exp(\beta_1 \widehat{x}+\beta_2 \widehat{y})\right)^{-1}.
\end{aligned}
\label{eq:local_field_parameterization}
\end{equation}
The logistic profile keeps $\widehat{b}$ strictly between zero and one, while the affine angle profile provides a compact representation of local directional variation. Together, these functions define the effective projected field
\begin{equation}
\widehat{\mathbf{b}}(\widehat{\mathbf{x}};\boldsymbol{\mu}_{\mathbf{B}})
=
\widehat{b}(\widehat{\mathbf{x}};\boldsymbol{\mu}_{\mathbf{B}})
\gamma(\widehat{\varphi}(\widehat{\mathbf{x}};\boldsymbol{\mu}_{\mathbf{B}})), 
\label{eq:local_vector_parameterization}
\end{equation}
where $\gamma: [0,2\pi] \to \mathbb{S}^1$ is the standard parametrization of the unit circle $\mathbb{S}^1$, i.e., 
\begin{equation}
    \gamma(\varphi) = (\cos(\varphi), \sin(\varphi)) \text{ for all } \varphi \in [0, 2\pi].
\end{equation}
The magnitude and orientation parameters therefore provide a low-dimensional surrogate for local variations in the magnetic geometry. Although this parameterization does not represent every possible magnetic equilibrium, it defines a family of smooth local tensor fields that can be sampled systematically.

For each parameter value, the corresponding diffusion tensor is defined by
\begin{equation}
\widehat{\mathbf{K}}(\boldsymbol{\mu}_{\mathbf{B}})
=
\kappa_{\perp}\mathbf{I}
+
(\kappa_{\parallel}-\kappa_{\perp})
\widehat{\mathbf{b}}(\cdot;\boldsymbol{\mu}_{\mathbf{B}})
\otimes
\widehat{\mathbf{b}}(\cdot;\boldsymbol{\mu}_{\mathbf{B}}),
\label{eq:local_tensor_parameterization}
\end{equation}
where the tensor is understood as a spatially varying matrix field on $\widehat{\omega}$. Thus, the parameterization induces a mapping from the finite-dimensional parameter space into a family of local diffusion operators.

\subsubsection{Snapshot Generation}
\label{subsubsec:snapshot_generation}

We discretize the parameter domain $\mathcal{P}_{\mathbf{B}}$ using a uniform Cartesian grid $\mathcal{G}_{\mathbf{B}}$. Let
$\left\{
\mathcal{H}_{\mathbf{B}}^{(m)}
\right\}_{m=1}^{N_m}$
denote the associated parameter-space cells. Each cell is characterized by its center 
$\left\{
\boldsymbol{\mu}_{\mathbf{B}}^{(m)}
\right\}_{m=1}^{N_m}
\subset\mathcal{P}_{\mathbf{B}}$, 
which will be used to construct a local reduced basis.

In addition to sampling the local tensor field, the offline procedure must account for the influence of boundary data on the local solution. This is important because a global solution restricted to a patch generally has nonzero traces on the patch boundary, even when the global problem has homogeneous Dirichlet conditions. We therefore generate local snapshots using a family of boundary excitations.
Let $T_{\partial\widehat{\omega}}:
\partial\widehat{\omega}\to[0,2\pi]$
be a boundary parameterization that traverses the perimeter of the reference patch. We interpret the coordinate $\theta=T_{\partial\widehat{\omega}}(\widehat{\mathbf{x}})$ as a periodic boundary coordinate. Its inverse is understood on the boundary, with the endpoints identified to represent the closed perimeter.
For a prescribed frequency cutoff $N_{\mathrm{freq}}$, we consider the trigonometric family
\begin{equation}
\begin{split}
q_{2j}(\theta) & = \cos(j\theta), \quad 0\leq j\leq N_{\mathrm{freq}},\\
q_{2j-1}(\theta) & = \sin(j\theta),\quad 1\leq j\leq N_{\mathrm{freq}}.
\end{split}
\label{eq:boundary_modes}
\end{equation}
This family contains $2N_{\mathrm{freq}}+1$ modes, including the constant mode. 
For $0 \leq j \leq 2 N_{\mathrm{freq}}$, the associated boundary function is defined by
\begin{equation}
\widehat{g}_j(\widehat{\mathbf{x}})
=
q_j\!\left(
T_{\partial\widehat{\omega}}(\widehat{\mathbf{x}})
\right),
\qquad \widehat{\mathbf{x}}\in\partial\widehat{\omega}.
\label{eq:boundary_trace_modes}
\end{equation}
The periodic representation ensures that the boundary values agree at the identified endpoints of the parameterization. In particular, the prescribed data are continuous around the perimeter when represented by the trigonometric functions in \eqref{eq:boundary_modes}.

Let $\widehat{V}_h$ denote a finite element space on the reference patch, and define its homogeneous-trace subspace by
\begin{equation}
\widehat{V}_{h,0}
=
\widehat{V}_h\cap H_0^1(\widehat{\omega}).
\end{equation}
For every parameter sample $\boldsymbol{\mu}_{\mathbf{B}}^{(m)}$, 
given the sample magnetic field 
$\widehat{\mathbf{b}}^{(m)} = \widehat{\mathbf{b}}(\bar{\boldsymbol{\mu}}_{\mathbf{B}}^{(m)})$ via \eqref{eq:local_vector_parameterization}, 
the sample diffusion tensor 
$\widehat{\mathbf{K}}^{(m)} = \widehat{\mathbf{K}}(\bar{\boldsymbol{\mu}}_{\mathbf{B}}^{(m)})$ is given by \eqref{eq:diffusion_tensor}.
For each boundary mode $j$, 
we compute a local snapshot $\widehat{w}_j^{(m)}$ with the prescribed boundary trace. Specifically, the snapshot $\widehat{w}_j^{(m)} \in \widehat{V}_h$ is the discrete harmonic extension 
\begin{equation}
\int_{\widehat{\omega}}
(\nabla \widehat{v})^\top
\widehat{\mathbf{K}}^{(m)}
\nabla\widehat{w}_j^{(m)}
\,\mathrm{d}\widehat{\mathbf{x}}
=0
\quad
\forall \widehat{v} \in \widehat{V}_{h,0},
\label{eq:snapshot_problem}
\end{equation}
satisfying the boundary condition
$\widehat{w}_j^{(m)}\big|_{\partial\widehat{\omega}}
=
\mathcal{J}_{\partial\widehat{\omega}}(\widehat{g}_j)$.
Here, $\mathcal{J}_{\partial\widehat{\omega}}$ denotes the boundary interpolation operator associated with the fine-scale discretization. The local snapshots are discrete harmonic extensions of the prescribed boundary data with respect to the parameterized anisotropic diffusion operator.

For each sampled magnetic configuration, this procedure produces a collection of local responses that captures how boundary variations are propagated through the anisotropic medium. Repeating the process over the parameter grid provides the snapshot data needed to construct reduced spaces that account for both variations in the local operator and variations in the local boundary trace.

\subsubsection{Spectral Reduction}
\label{subsubsec:spectral_reduction}

The collection of snapshots generated at all sampled parameter values can be substantially larger than the local reduced space required during the online phase. To compress these data while retaining representative local solution features, we perform a separate spectral reduction for each parameter-space cell $\mathcal{H}_{\mathbf{B}}^{(m)}$. We first construct the local snapshot space
\begin{equation}
\widehat{V}_{\mathrm{snap}}^{(m)}
=
\operatorname{span}
\left\{
\widehat{w}_j^{(m)}
\right\}_{0\leq j\leq 2N_{\mathrm{freq}}}.
\label{eq:snapshot_space}
\end{equation}
This space pools the local responses associated with the vertices surrounding the parameter-space cell. It therefore incorporates variations in both the local diffusion operator and the imposed boundary excitations.

We then define a local generalized eigenvalue problem using the diffusion operator evaluated at the cell center. Find eigenpairs 
$(\widehat{\lambda}_k^{(m)},\widehat{\phi}_k^{(m)})
\in\mathbb{R}\times \widehat{V}_{\mathrm{snap}}^{(m)}$
such that
\begin{equation}
\int_{\widehat{\omega}}
(\nabla \widehat{v})^\top
\widehat{\mathbf{K}}^{(m)}
\nabla\widehat{\phi}_k^{(m)}
\,\mathrm{d}\widehat{\mathbf{x}} =
\widehat{\lambda}_k^{(m)}
\int_{\widehat{\omega}}
\widehat{v} \widehat{\phi}_k^{(m)}
\,\mathrm{d}\widehat{\mathbf{x}}
\quad
\forall \widehat{v} \in \widehat{V}_{\mathrm{snap}}^{(m)}.
\label{eq:spectral_reduced}
\end{equation}
Equivalently, this problem identifies the modes of the snapshot space that diagonalize the local diffusion energy relative to the $L^2(\widehat{\omega})$ inner product. For a symmetric positive-definite tensor, the generalized eigenvalues are nonnegative, and positive eigenvalues correspond to nonzero diffusion energy.

We order the eigenvalues in ascending order
$0\leq\widehat{\lambda}_1^{(m)}
\leq\widehat{\lambda}_2^{(m)}
\leq\cdots \leq \widehat{\lambda}_{2N_{\mathrm{freq}} + 1}^{(m)}$
and the corresponding modes 
\begin{equation}
\widehat{\Phi}^{(m)}
=
\left\{
\widehat{\phi}_k^{(m)}
\right\}_{k=1}^{2N_{\mathrm{freq}} + 1}.
\label{eq:local_reduced_basis}
\end{equation}
Applying this procedure to all parameter-space cells yields the offline library
$\mathcal{L}
=
\left\{
\widehat{\Phi}^{(m)}
\right\}_{m=1}^{N_m}$,
which can be reused for global problems whose local tensor fields are adequately represented by the parameterized training family.

\subsection{PU-DDFEM solution scheme}
\label{subsec:online}

Let $\mathcal{S} = \{\omega_i\}_{i \in \mathcal{I}}$ be an overlapping open cover of $\Omega$, and
$\{\chi_i\}_{i \in \mathcal{I}} \subset H^1(\Omega)$ be the partition of unity
subordinate to $\mathcal{S}$ satisfying
\begin{equation}
\operatorname{supp}(\chi_i)\subseteq\overline{\omega_i},
\qquad
\sum_{i \in \mathcal{I}}\chi_i(\mathbf{x})=1
\quad\text{for }\mathbf{x}\in\Omega.
\label{eq:pou_properties}
\end{equation}
We select a subset of patches $\{\omega_i\}_{i \in \mathcal{I}_{\text{DOF}}}$
indexed by $\mathcal{I}_{\text{DOF}}$, for each
we identify a suitable local model from the precomputed library which best approximates the local magnetic geometry in $\omega_i$, 
conformalize the basis functions on the fine background mesh $\mathcal{T}_h$ through the partition of unity, 
and assemble the global reduced space $V_{H,0} \subset V_{h,0}$ for solution approximation.
We remark that constructing both the partition of unity $\{\chi_i\}_{i \in \mathcal{I}}$ and the active index set $\mathcal{I}_{\text{DOF}}$ requires careful coordination to simultaneously preserve the partition of unity property and enforce the required essential boundary conditions.

\subsubsection{Local model identification}
\label{subsubsec:model_identification}

For a given patch $\omega_i$, where $i \in \mathcal{I}_{\text{DOF}}$, 
using the geometric mapping 
$T_i:\omega_i\to\widehat{\omega}$, 
the local material tensor field can be expressed in reference coordinates as
\begin{equation}
\widehat{\mathbf{K}}_i(\widehat{\mathbf{x}})
=
\mathbf{K}(T_i^{-1}(\widehat{\mathbf{x}})),
\qquad
\widehat{\mathbf{x}}\in\widehat{\omega}.
\label{eq:reference_tensor}
\end{equation}
We identify the library cell whose representative tensor most closely matches this mapped field. Specifically, we select
\begin{equation}
m_i^*
=
\underset{1\leq m\leq N_m}{\operatorname{argmin}}
\left\|
\widehat{\mathbf{K}}_i
-
\widehat{\mathbf{K}}(\bar{\boldsymbol{\mu}}_{\mathbf{B}}^{(m)})
\right\|_{L^2(\widehat{\omega})},
\label{eq:local_matching}
\end{equation}
where the tensor-field norm is defined by
\begin{equation}
\|\widehat{\mathbf{E}}\|_{L^2(\widehat{\omega})}
=
\left(
\int_{\widehat{\omega}}
\|\widehat{\mathbf{E}}(x)\|_F^2
\,\mathrm{d}\widehat{\mathbf{x}}
\right)^{1/2},
\end{equation}
and $\|\cdot\|_F$ denotes the Frobenius norm of a matrix.
This selection procedure provides a direct way to identify a precomputed local model without solving an additional local spectral problem. Its accuracy relies on the representativeness of the training parameter space and on the extent to which proximity in the tensor norm predicts similarity of the associated local solution spaces. In particular, the matching criterion is a model-selection heuristic rather than an a priori error estimator for the resulting global solution.
Once the library is identified, we can choose to retain the first $r_i$ modes, where $1 \leq r_i \leq 2N_{\mathrm{freq}} + 1$. 
These modes span the selected low-energy portion of the local snapshot space, and the quality of this truncation depends on the local tensor, the snapshot ensemble, and the selected dimension.

\subsubsection{Global assembly}
\label{subsubsec:global_assembly}

After identifying the library index $m_i^*$ and choosing the dimension $r_i$ for each patch, we transfer the selected reference basis functions to their corresponding physical subdomains. 
For a selected local mode $\widehat{\phi}_k^{(m_i^*)}$, define the corresponding physical basis function by 
$\psi_{i,k} = \mathcal{J}_{\Omega}(\widetilde{\psi}_{i,k}) \in V_{h,0}$,
where  
\begin{equation}
\widetilde{\psi}_{i,k}(\mathbf{x})
=
\begin{cases}
\chi_i(\mathbf{x})\,
\widehat{\phi}_k^{(m_i^*)}(T_i(\mathbf{x})),
&\mathbf{x}\in\omega_i,\\
0,
&\mathbf{x}\in\Omega \setminus \omega_i,
\end{cases}
\label{eq:global_basis}
\end{equation}
for $1\leq k\leq r_i$. 
We consequently define the global reduced space as
\begin{equation}
V_{H,0}
=
\operatorname{span}
\left\{
\psi_{i,k}:
1\leq k\leq r_i, \;
i \in \mathcal{I}_{\text{DOF}}
\right\} \subset V_{h,0}.
\label{eq:global_reduced_space}
\end{equation}
The dimension of the assembled PU-DDFE space is
\begin{equation}
N_H = \dim(V_{H,0})
= \sum_{i \in \mathcal{I}_{\text{DOF}}} r_i.
\label{eq:reduced_dimension}
\end{equation}

\subsubsection{PU-DDFEM solution}
\label{subsubsec:PU-DDFEM-sol}
The PU-DDFEM solution $u_H = u_{h,g} + u_{H,0} \in V_h$ 
is then obtained through a reduced Galerkin formulation: find $u_{H,0} \in V_{H,0}$ such that 
\begin{equation}
a(u_{H,0},v)
=
\ell(v) - a(u_{h,g}, v)
\quad
\forall v\in V_{H,0}.
\label{eq:PU-DDFEM-sol}
\end{equation}
By re-indexing the double-indexed basis functions into a single index, the global reduced basis can be written as $\{\psi_j\}_{j=1}^{N_H}$. Since $V_{H,0} \subset V_{h,0}$, each reduced basis function $\psi_j$ can be expressed as a linear combination of the high-fidelity shape basis functions $\{v_k\}_{k=1}^{N_h}$ via
\begin{equation}
\psi_j = \sum_{k=1}^{N_h} \Psi_{kj} v_k, \quad j = 1, \dots, N_H,
\end{equation}
where $\boldsymbol{\Psi} = (\Psi_{kj})_{1 \leq k \leq N_h, 1 \leq j \leq N_H} \in \mathbb{R}^{N_h \times N_H}$ denotes the transformation matrix containing these expansion coefficients. Expressing the reduced solution as a linear combination of the reduced basis, \eqref{eq:PU-DDFEM-sol} can be expressed as reduced linear system via Galerkin projection
\begin{equation}
\mathbf{A}_H \mathbf{U}_H = \mathbf{F}_H,
\label{eq:reduced_linear}
\end{equation}
where $\mathbf{A}_H = \boldsymbol{\Psi}^T \mathbf{A}_h \boldsymbol{\Psi} \in \mathbb{R}^{N_H \times N_H}$ is the reduced stiffness matrix, $\mathbf{F}_H = \boldsymbol{\Psi}^T \mathbf{F}_h \in \mathbb{R}^{N_H}$ is the reduced load vector, and $\mathbf{U}_H \in \mathbb{R}^{N_H}$ represents the vector of reduced coordinates. This reduced Galerkin problem retains the variational structure of the weak form \eqref{eq:weak_form}, while seeking an accurate approximation in a substantially smaller space.
Thanks to the block sparse structure, symmetry, and positive definiteness of the matrix $\mathbf{A}_H$, the linear system \eqref{eq:reduced_linear} is solved using the sparse LDL factorization.
The computational advantage arises when $N_H$ is substantially smaller than $N_h$ and using the PU-DDFE space for solving multiple source functions. The precomputed library avoids repeating local snapshot generation and spectral reduction for every global problem, while the reduced Galerkin system decreases the dimension of the online solve. The factorization is pre-computed for a given magnetic flux configuration and re-usable for different source functions, following the common practice of offline-online splitting in reduced order modeling. The total computational benefit nevertheless depends on the cost of local model identification, basis transfer, global assembly, and reduced-system solution, as well as the accuracy achieved for the selected local basis dimensions.
\section{Numerical Examples}
\label{sec:numerical}

To evaluate the accuracy and robustness of the proposed PU-DDFEM, we present a series of numerical benchmarks motivated by heat and particle transport in magnetically confined fusion plasmas. Specifically, we investigate three iconic magnetic field topologies: the single-null, double-null, and magnetic island configurations.
All the simulations in this section use Dane in Livermore Computing Center\footnote{High
performance computing at LLNL, https://hpc.llnl.gov/hardware/compute-platforms/dane},
on Intel Sapphire Rapids CPUs with 256 GB memory/node and peak PFLOPS of 10.7.

To quantify the solution accuracy, we define the relative $L^2$ error between the PU-DDFEM solution and the reference solution
\begin{equation}
\varepsilon_{L^2} = \dfrac{\| u_H - u_h \|_{L^2(\Omega)}}{\| u_h \|_{L^2(\Omega)}}.
\label{eq:L2_error}
\end{equation}
We also report the relative wall clock time $\tau$, measured by the ratio of the wall clock time of $u_H$ to that of $u_h$.

\subsection{Benchmark problems and FEM solution}
\label{subsec:benchmarks}

The anisotropic thermal diffusion tensor $\mathbf{K}$ is prescribed analytically via scalar magnetic flux function $A_z$, 
the perpendicular diffusion coefficient $\kappa_{\perp} = 1$,
and the parallel diffusion coefficient $\kappa_{\parallel} \in \{10^2, 10^4, 10^6\}$. 
Heat transport occurs predominantly along the magnetic field lines while being strongly suppressed across them. 
Consequently, the local topology of the magnetic field dictates the heat flow pathways and governs the formation of sharp spatial temperature gradients. 
In regions featuring closed magnetic field lines surrounded by a separatrix, heat introduced via volumetric sources or non-homogeneous boundary conditions remains largely trapped within the core, establishing internal boundary layers of characteristic width $\mathcal{O}(\kappa_{\parallel}^{-1/2})$ across the separatrix interface. These localized internal layers present severe challenges for standard numerical methods on unaligned grids, making them ideal benchmarks for our proposed approach. 

In all benchmark test cases, the computation domain is $\Omega = (0,1)^2$. 
Throughout this work, we consider rectangular partitions. 
The domain $\Omega$ is partitioned into a uniform mesh $\mathcal{T}_h$ of square elements with length $h$. 
The mesh parameter $h$ is chosen depending on the parallel diffusion coefficient $\kappa_{\parallel}$, 
namely $h = 2^{-(7+\log_{10} \sqrt{\kappa_{\parallel}})}$ for $\kappa_{\parallel} \in \{10^2, 10^4, 10^6\}$, 
to ensure the fine-scale details are sufficiently captured in the corresponding anisotropy level.
The finite element space $V_h$ consisting of $H^1$-conforming bilinear elements is defined as:
$$V_h = \{ v_h \in H^1(\Omega) : v_h\vert_K \in \mathcal{Q}_1(K)\, \forall K \in \mathcal{T}_h \}$$
where $\mathcal{Q}_1(K)$ denotes the space of bilinear polynomials on the element $K$. 

\subsubsection{Single-null configuration}
\label{subsubsec:single_null}

\begin{figure}[htbp!]
    \centering
    \includegraphics[width=0.48\linewidth]{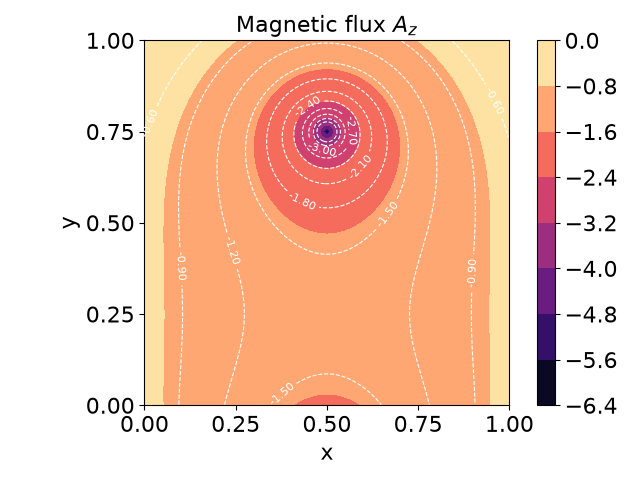}
    \includegraphics[width=0.48\linewidth]{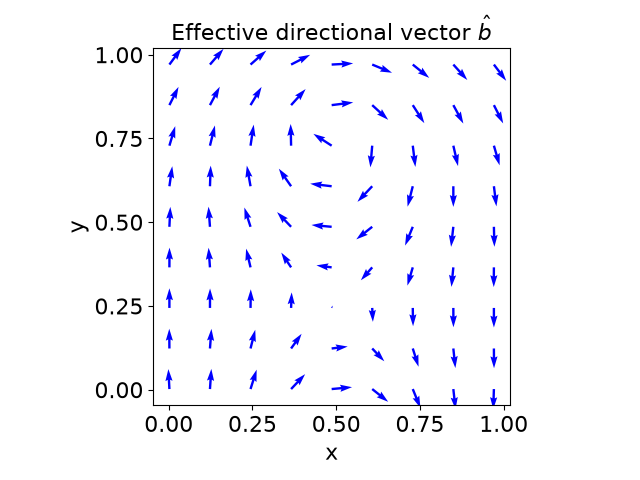}
    \caption{Magnetic flux function $A_z$ (left) and the resultant effective directional vector $\widehat{\mathbf{b}}$ (right) in the single-null configuration.}
    \label{fig:single_null_magnetic}
\end{figure}

\begin{figure}[htbp!]
    \centering
    \includegraphics[width=0.48\linewidth]{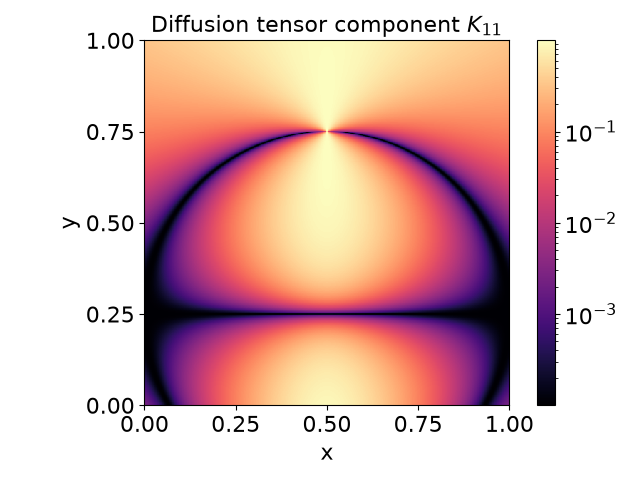}
    \includegraphics[width=0.48\linewidth]{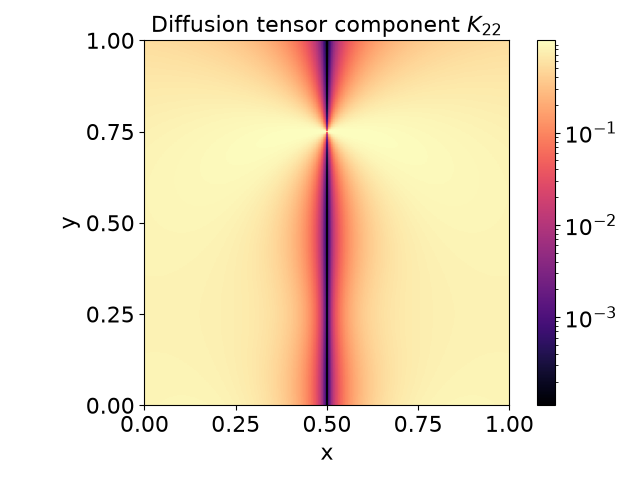}
    \caption{Components $K_{11}$ (left) and $K_{22}$ (right) of diffusion tensor $\mathbf{K}$ in the single-null configuration with anisotropic ratio $\kappa_{\parallel} = 10^4$.}
    \label{fig:single_null_diffusion}
\end{figure}

The first benchmark considers a single-null magnetic field topology, representative of standard diverted tokamak configurations. This topology contains a single magnetic poloidal null point that partitions the domain into closed flux surfaces in the core and open field lines in the scrape-off layer intersecting the domain boundaries. In the computational domain $\Omega = (0,1)^2$, the magnetic flux function $A_z$ is generated via a two-wire current system given by: for all $\mathbf{x} = (x,y) \in \Omega$,
\begin{equation}
A_z(\mathbf{x}) = \dfrac{1}{2} \left(\log\left((x-0.5)^2+(y-0.75)^2\right) + \log\left((x-0.5)^2+(y+0.25)^2\right)\right).
\end{equation}
The resulting magnetic field features a single poloidal null point at $\mathbf{x}_s = (0.5, 0.25)$. 
Figure~\ref{fig:single_null_magnetic} illustrates the contours of the magnetic flux function $A_z$ alongside the quiver plot of resultant effective directional vector $\widehat{\mathbf{b}}$ across the computational domain $\Omega$. Figure~\ref{fig:single_null_diffusion} illustrates the components $K_{11}$ and $K_{22}$ of the diffusion tensor $\mathbf{K}$ with anisotropic ratio $\kappa_{\parallel} = 10^4$. 

To model core heat generation, a localized volumetric heat source $f$ is placed inside the closed separatrix region, represented by the isotropic Gaussian profile: for all $\mathbf{x} = (x,y) \in \Omega$,
\begin{equation}
\label{eq:gaussian_source}
f(\mathbf{x}) = \exp\left( -\frac{1}{2} \left[ (8(x-0.5))^2 + (8(y-0.5))^2 \right] \right),
\end{equation}
which is centered at $(0.5, 0.5)$ and decays rapidly toward the domain boundaries. 
Homogeneous Dirichlet boundary conditions $u = 0$ are imposed on $\partial\Omega$.
The affine space in the weak formulation \eqref{eq:weak_form} is reduced to a linear space 
\begin{equation}
u_g = 0, \quad V_0 = \left\{ v \in H^1(\Omega) : v = 0\right\}, 
\label{eq:null_problem_space}
\end{equation}
where the condition is understood in the sense of the trace operator. 
Similarly, the affine space in the finite element formulation \eqref{eq:fem} is reduced to a linear space with $u_{h,g} = 0$.

\begin{figure}[htbp!]
    \centering
    \includegraphics[width=0.32\linewidth]{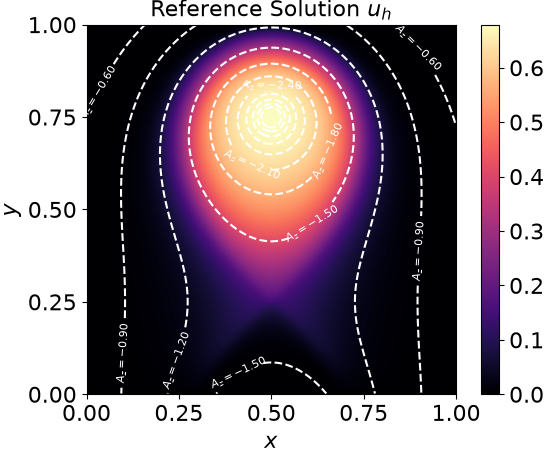}
    \includegraphics[width=0.32\linewidth]{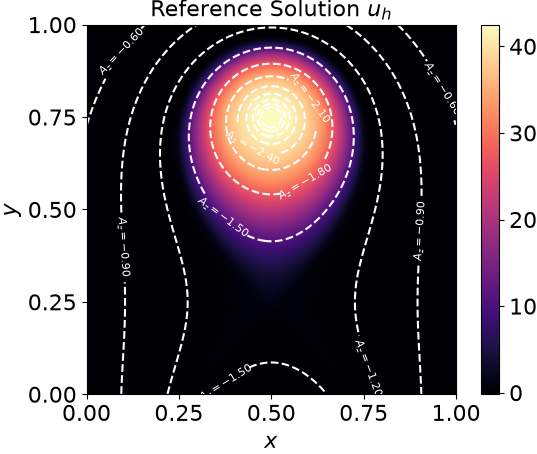}
    \includegraphics[width=0.32\linewidth]{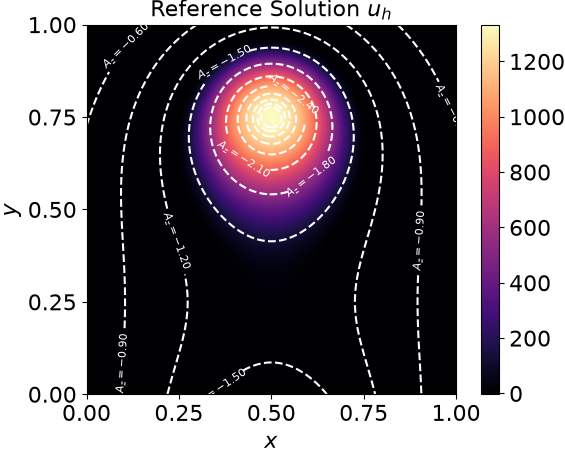}
    \caption{Reference FEM solution $u_h$ for the single-null problem with anisotropic ratio $\kappa_{\parallel} = 10^2$ (left), $\kappa_{\parallel} = 10^4$ (middle), and $\kappa_{\parallel} = 10^6$ (right).}
    \label{fig:single_null_fem_solution}
\end{figure}

The FEM solution $u_h$ in \eqref{eq:fem} for the single-null problem with anisotropic ratio $\kappa_{\parallel} \in \{10^2, 10^4, 10^6\}$ are shown in Figure~\ref{fig:single_null_fem_solution}. 
The contours of the magnetic flux function $A_z$ are plotted to illustrate the separatrix. 
Since parallel diffusion along the closed field lines dominates perpendicular transport, heat rapidly equilibrates along flux contours while forming steep internal layers across the separatrix interface.

\subsubsection{Double-null configuration}
\label{subsubsec:double_null}

\begin{figure}[htbp!]
    \centering
    \includegraphics[width=0.48\linewidth]{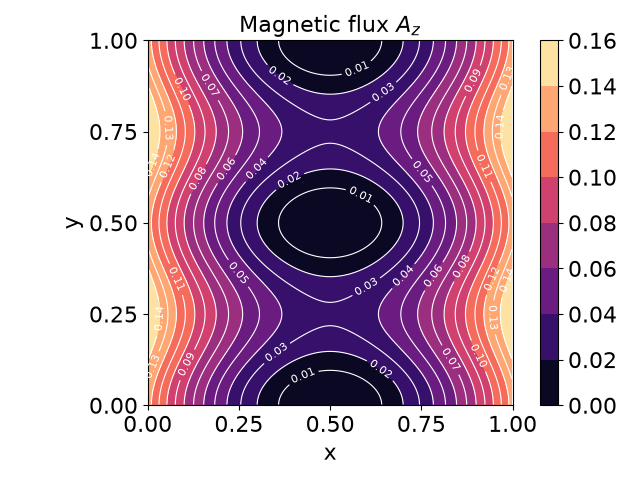}
    \includegraphics[width=0.48\linewidth]{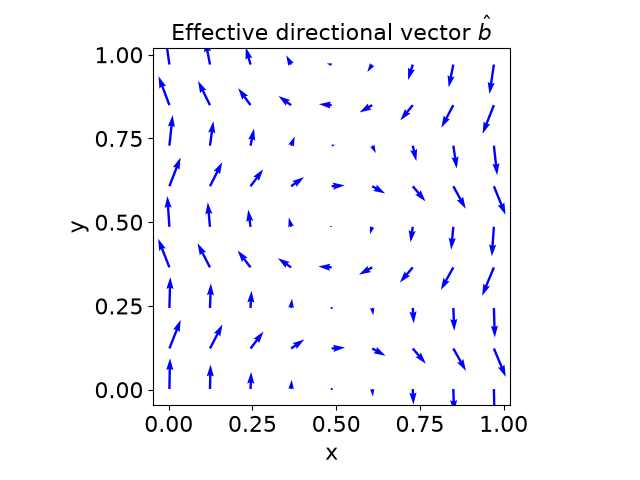}
    \caption{Magnetic flux function $A_z$ (left) and the resultant effective directional vector $\widehat{\mathbf{b}}$ (right) in the double-null configuration.}
    \label{fig:double_null_magnetic}
\end{figure}

\begin{figure}[htbp!]
    \centering
    \includegraphics[width=0.48\linewidth]{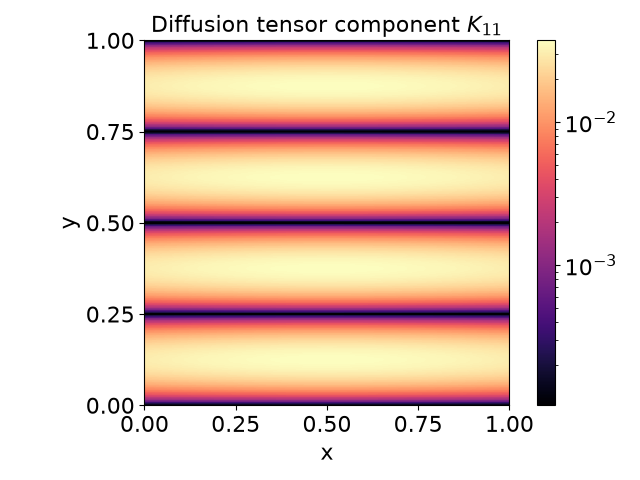}
    \includegraphics[width=0.48\linewidth]{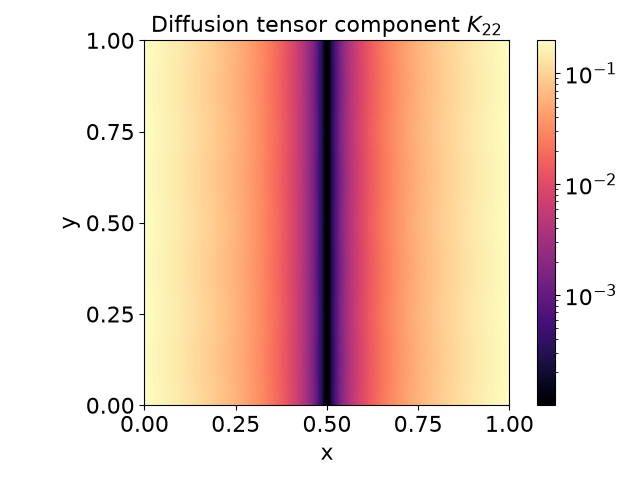}
    \caption{Components $K_{11}$ (left) and $K_{22}$ (right) of diffusion tensor $\mathbf{K}$ in the double-null configuration with anisotropic ratio $\kappa_{\parallel} = 10^4$.}
    \label{fig:double_null_diffusion}
\end{figure}

The second benchmark considers a double-null magnetic topology, which features two symmetric magnetic null points. This configuration provides a more challenging, non-linear field line curvature distribution and complex magnetic shear. The scalar flux function $A_z$ is defined by: for all $\mathbf{x} = (x,y) \in \Omega = (0,1)^2$ as:
\begin{equation}
A_z(\mathbf{x}) = \frac{1}{2}(x-0.5)^2 + \frac{1}{2} \left( \frac{1}{4}\sin\left(2\pi(y-0.5)\right) \right)^2.
\end{equation}
This formulation yields two magnetic nulls located at $\mathbf{x}_s^{\pm} = (0.5, 0.5 \pm 0.25)$, connected by intersecting separatrices that separate the inner core, top/bottom divertor channels, and outer edge regions. Figure~\ref{fig:double_null_magnetic} depicts the magnetic flux function $A_z$ and the effective directional vector $\widehat{\mathbf{b}}$ for this configuration.
Figure~\ref{fig:double_null_diffusion} depicts the diagonal components $K_{11}$ and $K_{22}$ of the diffusion tensor $\mathbf{K}$ with anisotropic ratio $\kappa_{\parallel} = 10^4$. 

As in the single-null problem, we prescribe the same Gaussian volumetric heat source $f$ given by \eqref{eq:gaussian_source} centered at $(0.5, 0.5)$, combined with homogeneous Dirichlet boundary conditions $u = 0$ along $\partial\Omega$. 

\begin{figure}[htbp!]
    \centering
    \includegraphics[width=0.32\linewidth]{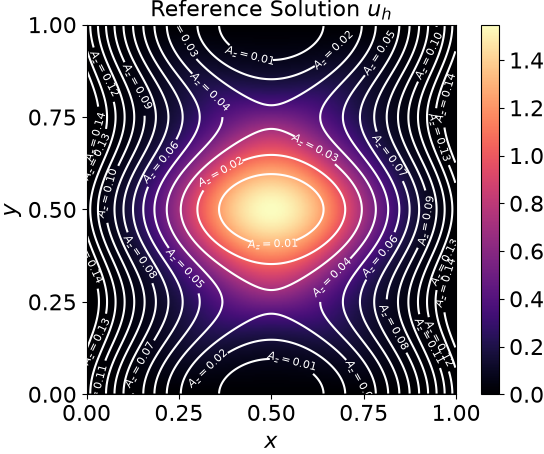}
    \includegraphics[width=0.32\linewidth]{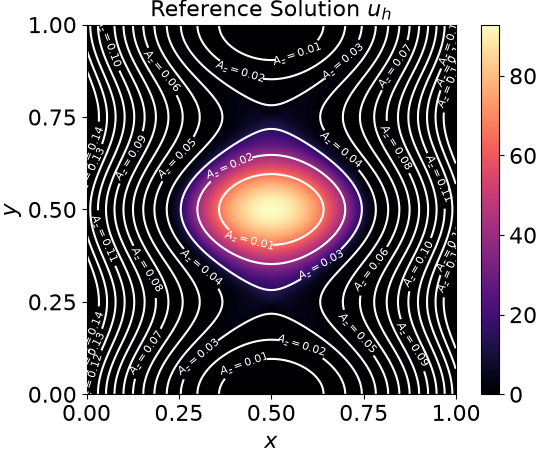}
    \includegraphics[width=0.32\linewidth]{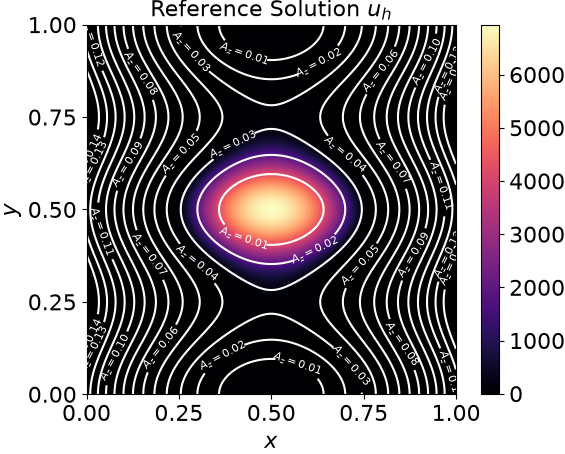}
    \caption{Reference FEM solution $u_h$ for the double-null problem with anisotropic ratio $\kappa_{\parallel} = 10^2$ (left), $\kappa_{\parallel} = 10^4$ (middle), and $\kappa_{\parallel} = 10^6$ (right).}
    \label{fig:double_null_fem_solution}
\end{figure}

Figure~\ref{fig:double_null_fem_solution} shows the FEM solution in \eqref{eq:fem} for the double-null problem with anisotropic ratio $\kappa_{\parallel} \in \{10^2, 10^4, 10^6\}$. 
The presence of dual null points creates complex stagnation regions where $\mathbf{B}_p$ vanishes, intensifying the numerical cross-wind pollution on unaligned grids and testing the capability of our pre-trained local basis library to capture intersecting internal layers.

\subsubsection{Magnetic island configuration}
\label{subsubsec:magnetic_island}

\begin{figure}[htbp!]
    \centering
    \includegraphics[width=0.48\linewidth]{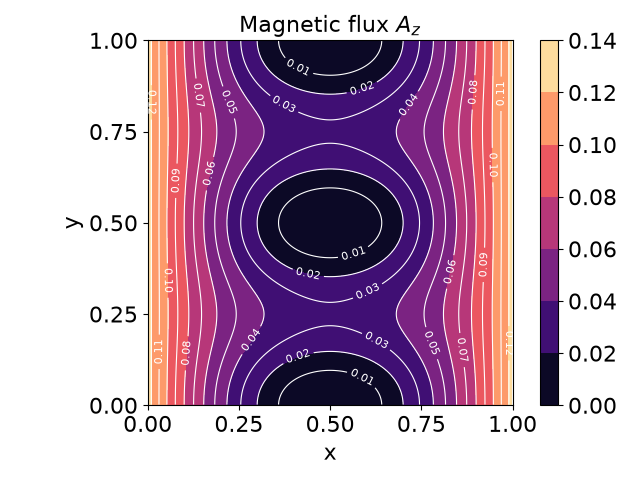}
    \includegraphics[width=0.48\linewidth]{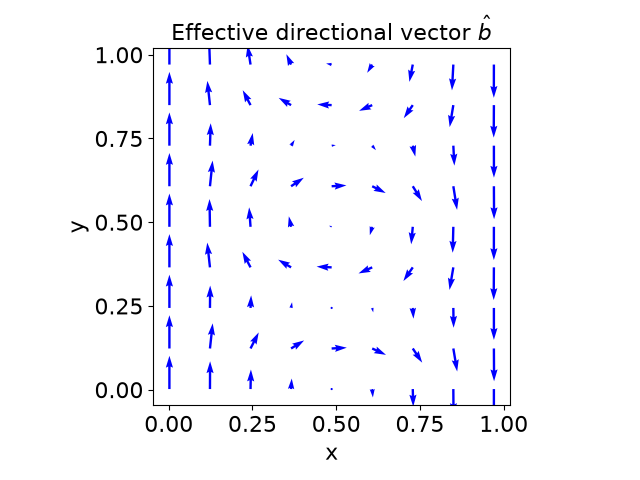}
    \caption{Magnetic flux function $A_z$ (left) and the resultant effective directional vector $\widehat{\mathbf{b}}$ (right) in the double-null configuration.}
    \label{fig:magnetic_island_magnetic}
\end{figure}

\begin{figure}[htbp!]
    \centering
    \includegraphics[width=0.48\linewidth]{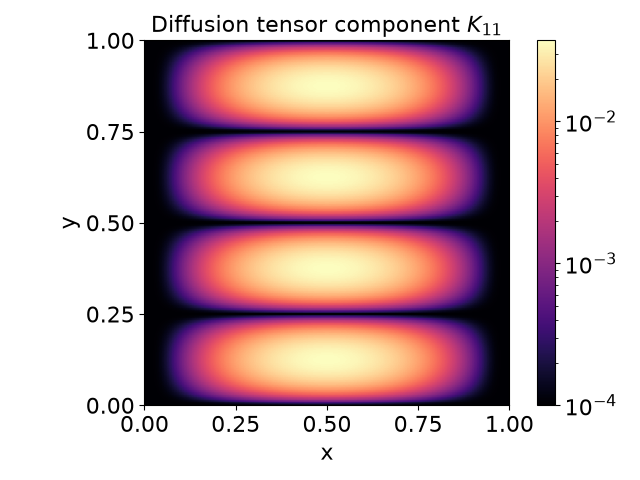}
    \includegraphics[width=0.48\linewidth]{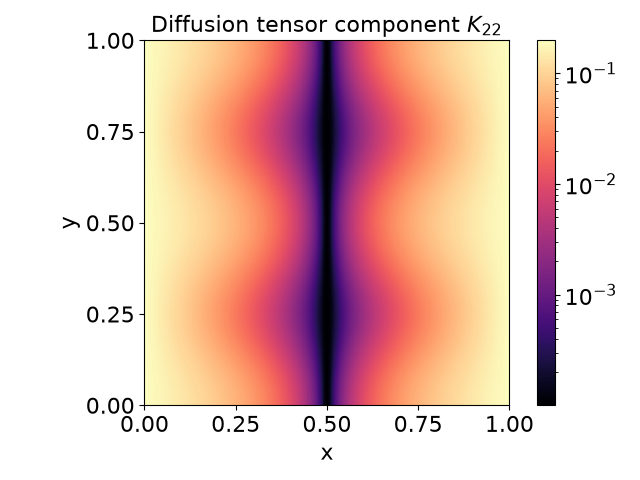}
    \caption{Components $K_{11}$ (left) and $K_{22}$ (right) of diffusion tensor $\mathbf{K}$ in the double-null configuration with anisotropic ratio $\kappa_{\parallel} = 10^4$.}
    \label{fig:magnetic_island_diffusion}
\end{figure}

The third benchmark investigates a magnetic island topology, which models how resonant magnetic perturbations disrupt nested magnetic flux surfaces, inducing field line reconnection and forming closed island chains within the plasma volume. The flux function $A_z$ incorporates a spatially tapered sinusoidal perturbation added to a background equilibrium field, defined on $\Omega = (0,1)^2$ by: for all $\mathbf{x} = (x,y) \in \Omega$,
\begin{equation}
A_z(\mathbf{x}) = \frac{1}{2}(x-0.5)^2 + \frac{1}{32}\left( \left(1 - 4(x-0.5)^2\right) \sin\left(2\pi(y-0.5)\right) \right)^2.
\label{eq:mag_island_flux}
\end{equation}
The quadratic term $\frac{1}{2}(x-0.5)^2$ represents the background sheared magnetic field, while the factor $(1 - (x-0.5)^2 / 0.25)$ enforces the perturbation to vanish smoothly at the radial boundaries $x=0$ and $x=1$. This perturbation opens a magnetic island centered at $(0.5, 0.5)$, bounded by a local separatrix. Figure~\ref{fig:magnetic_island_magnetic} illustrates the flux function $A_z$ and the resulting magnetic vector $\widehat{\mathbf{b}}$.
Figure~\ref{fig:magnetic_island_diffusion} demonstrates the diffusion coefficients $K_{11}$ and $K_{22}$ with anisotropic ratio $\kappa_{\parallel} = 10^4$.

Unlike the single- and double-null benchmarks driven by internal volumetric sources, the magnetic island problem is driven by a cross-domain thermal gradient with inhomogeneous Dirichlet boundary conditions on the radial boundaries and periodic boundary conditions in the poloidal direction, i.e., 
\begin{equation}
\begin{split}
u(0, y) = 0, \; u(1, y) = 1, & \quad \forall y \in [0,1], \\
u(x, 0) = u(x, 1), \; \frac{\partial u}{\partial y}(x, 0) = \frac{\partial u}{\partial y}(x, 1), & \quad \forall x \in [0,1].
\end{split}
\end{equation}
The boundary condition maintains an effective heat source along the hot boundary $x=0$, and thermal energy must transport across the domain to the cold boundary $x=1$.
The affine space for weak formulation in \eqref{eq:weak_form} is given by 
\begin{equation}
u_g(\mathbf{x}) = 1-x, \quad 
V_0 = \left\{ v \in H^1(\Omega) \;\middle|\;
\begin{aligned}
&v(0, y) = 0, \; v(1, y) = 0 && \forall y \in [0,1] \\
&v(x, 0) = v(x, 1) && \forall x \in [0,1]
\end{aligned}
\right\}.
\label{eq:magnetic_island_space}
\end{equation}
Since $u_g \in V_h$, the affine space in the finite element formulation \eqref{eq:fem} has the origin $u_{h,g} = u_g$.

\begin{figure}[htbp!]
    \centering
    \includegraphics[width=0.32\linewidth]{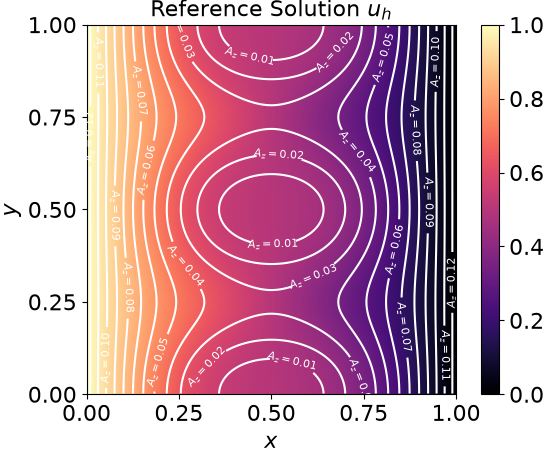}
    \includegraphics[width=0.32\linewidth]{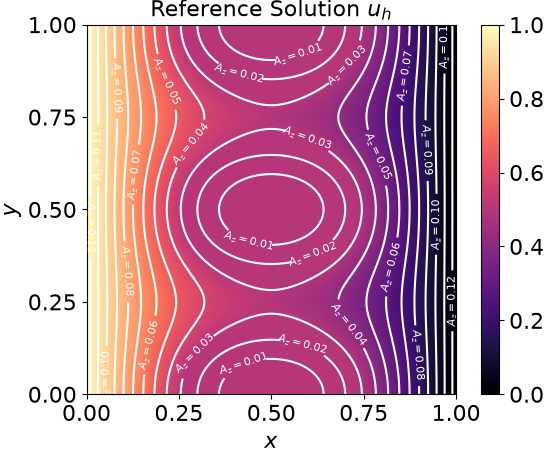}
    \includegraphics[width=0.32\linewidth]{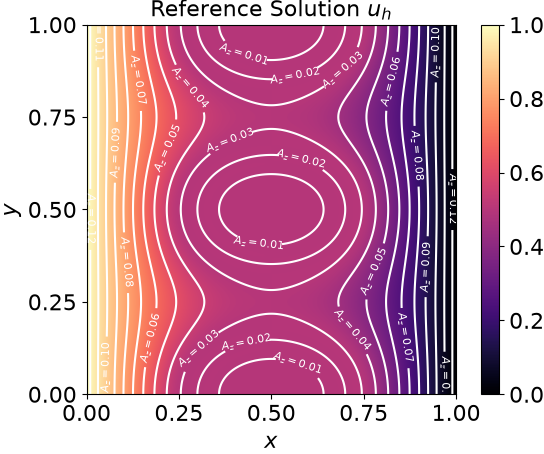}
    \caption{Reference FEM solution $u_h$ for the magnetic island problem with anisotropic ratio $\kappa_{\parallel} = 10^2$ (left), $\kappa_{\parallel} = 10^4$ (middle), and $\kappa_{\parallel} = 10^6$ (right).}
    \label{fig:magnetic_island_fem_solution}
\end{figure}

Figure~\ref{fig:magnetic_island_fem_solution} depicts the FEM solution in \eqref{eq:fem} for the double-null problem with anisotropic ratio $\kappa_{\parallel} \in \{10^2, 10^4, 10^6\}$. 
With high parallel conductivity $\kappa_{\parallel}$, heat flattens rapidly within the closed flux surfaces of the island interior, creating a constant-temperature plateau inside the island core bounded by steep localized boundary layers along the island separatrix.

\subsection{PU-DDFEM construction}

\begin{figure}
    \centering
    \includegraphics[width=0.32\linewidth]{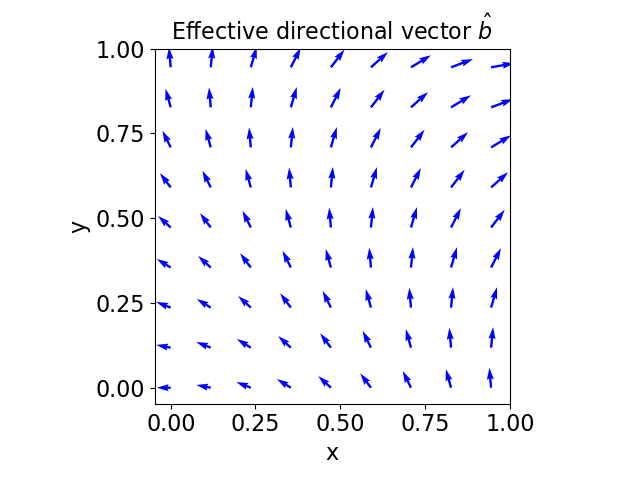}
    \includegraphics[width=0.32\linewidth]{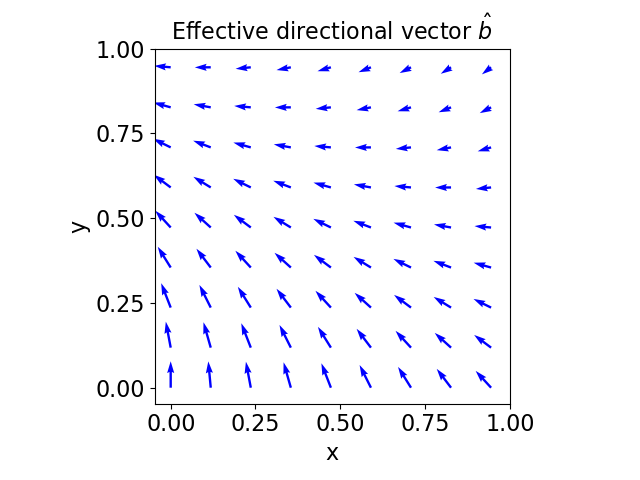}
    \includegraphics[width=0.32\linewidth]{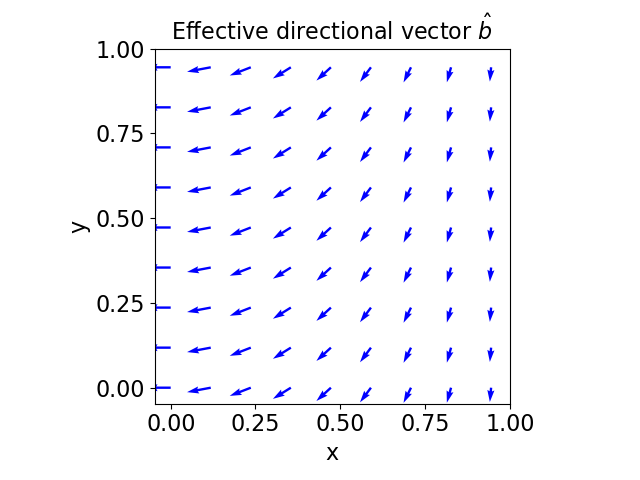}\\
    \includegraphics[width=0.32\linewidth]{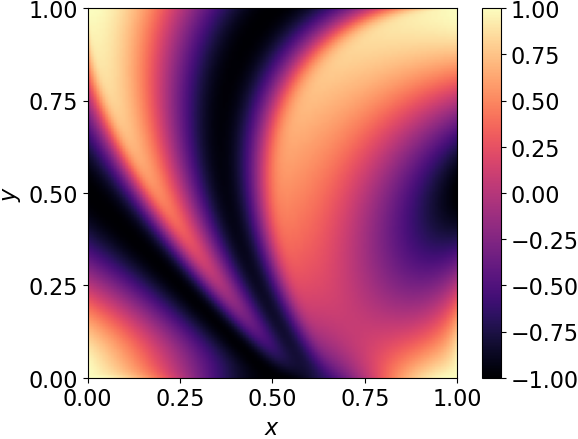}
    \includegraphics[width=0.32\linewidth]{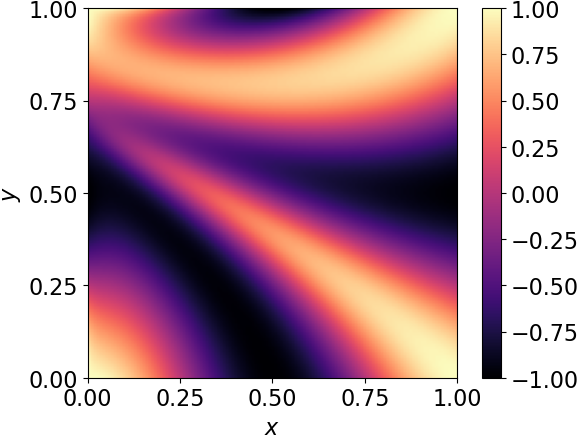}
    \includegraphics[width=0.32\linewidth]{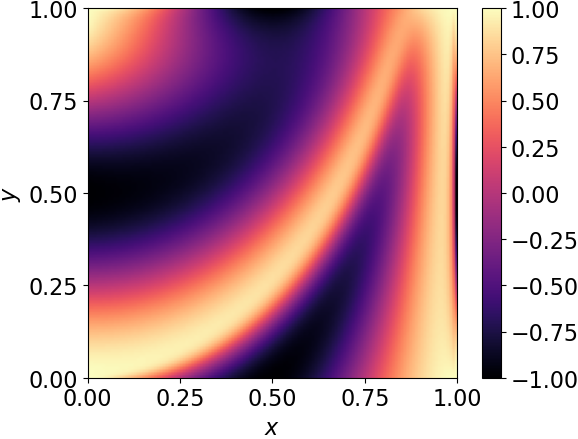}\\
    \includegraphics[width=0.32\linewidth]{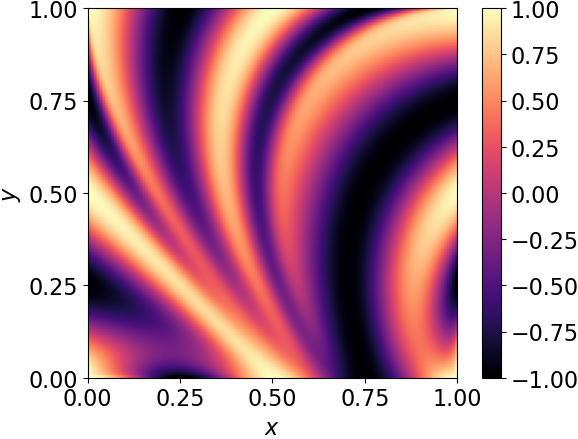}
    \includegraphics[width=0.32\linewidth]{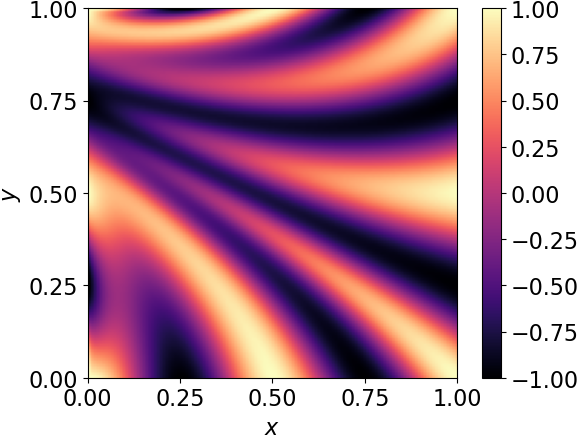}
    \includegraphics[width=0.32\linewidth]{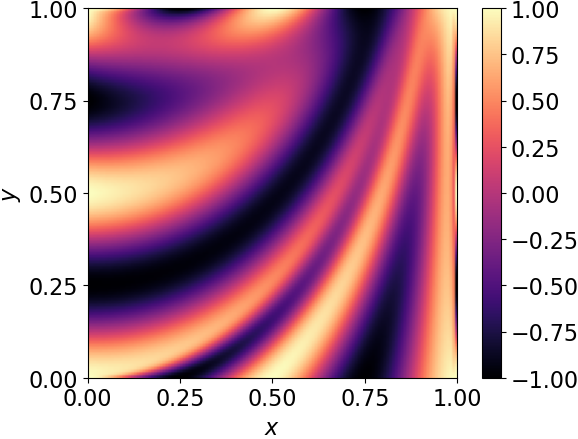}\\
    \includegraphics[width=0.32\linewidth]{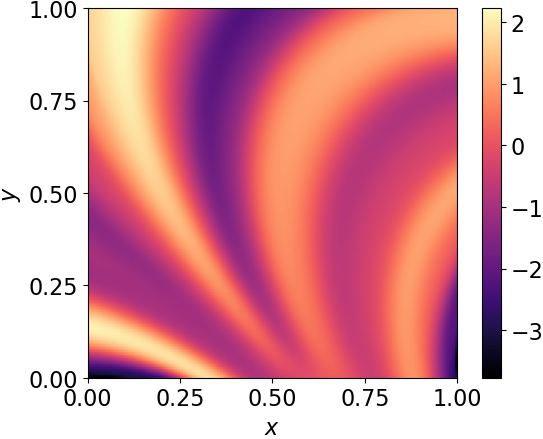}
    \includegraphics[width=0.32\linewidth]{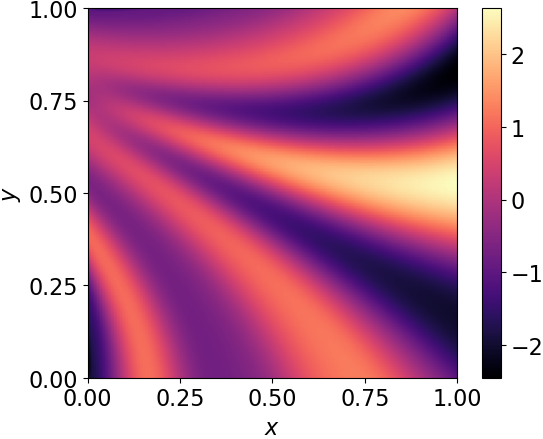}
    \includegraphics[width=0.32\linewidth]{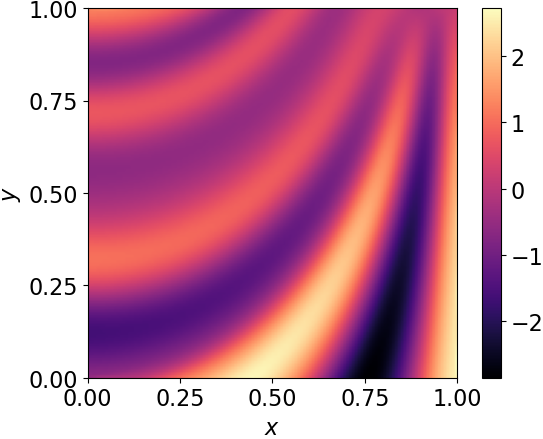}\\
    \includegraphics[width=0.32\linewidth]{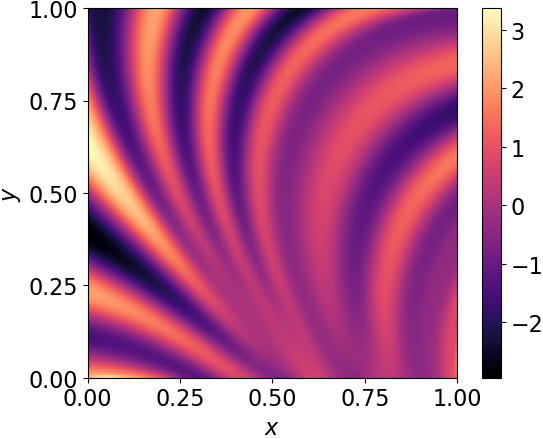}
    \includegraphics[width=0.32\linewidth]{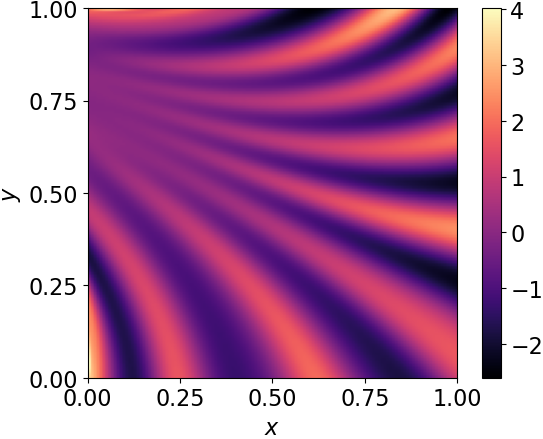}
    \includegraphics[width=0.32\linewidth]{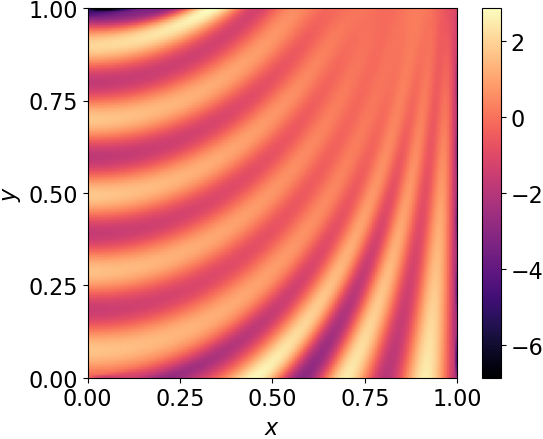}\\
    \caption{Representative samples of parameterized magnetic flux vector fields $\widehat{\mathbf{b}}^{(m)}$ (1st row), the corresponding snapshots $\widehat{w}_{8}^{(m)}, \widehat{w}_{16}^{(m)}$ (2nd and 3rd row), and the corresponding basis $\widehat{\phi}_{8}^{(m)}, \widehat{\phi}_{16}^{(m)}$ (4th and 5th row) on the reference patch $\widehat{\omega}$.}
    \label{fig:sample}
\end{figure}

We consider the reference patch $\widehat{\omega} = (0,1)^2$ and 
the parametric domain  
$$\mathcal{P}_{\mathbf{B}} = [-\pi, \pi] \times \left[-\frac{\pi}{2}, \frac{\pi}{2}\right] \times \left[-\frac{\pi}{2}, \frac{\pi}{2}\right] \times [-1, 1] \times [-1, 1].$$
in which we take combinations of  
\begin{equation}
    \begin{split}
    \varphi_0 & \in \left\{-\pi, -\dfrac{\pi}{2}, 0, \dfrac{\pi}{2}, \pi\right\}, \\
    \alpha_1, \alpha_2 & \in \left\{-\dfrac{\pi}{2}, -\dfrac{\pi}{4}, 0, \dfrac{\pi}{4}, \dfrac{\pi}{2}\right\}, \\
    \beta_1, \beta_2 & \in \left\{-1, -\dfrac{1}{2}, 0, \dfrac{1}{2}, 1\right\}, \\
    \end{split}
\end{equation}
for constructing the local magnetic fields via \eqref{eq:local_field_parameterization}, 
computing the harmonic extensions of the boundary modes as stated in \eqref{eq:snapshot_problem},  
and performing spectral reduction as described in \eqref{eq:spectral_reduced}. 
Figure~\ref{fig:sample} illustrates 3 representative samples of 
the parameterized magnetic field $\widehat{\mathbf{b}}^{(m)}$, 
the corresponding snapshots 
$\widehat{w}_{8}^{(m)}, \widehat{w}_{16}^{(m)}$
defined in \eqref{eq:snapshot_problem} with $N_\text{freq} = 20$, and
the corresponding basis 
$\widehat{\phi}_{8}^{(m)}, \widehat{\phi}_{16}^{(m)}$
defined in \eqref{eq:spectral_reduced}, for $m = m_1, m_2, m_3$, where 
\begin{equation}
    \begin{split}
    \boldsymbol{\mu}_{\mathbf{B}}^{(m_1)} & = \left(-\pi, -\dfrac{\pi}{2}, -\dfrac{\pi}{2}, -1, -1 \right), \\
    \boldsymbol{\mu}_{\mathbf{B}}^{(m_2)} & = \left(\dfrac{\pi}{2}, \dfrac{\pi}{4}, \dfrac{\pi}{2}, \dfrac{1}{2}, 1\right), \\
    \boldsymbol{\mu}_{\mathbf{B}}^{(m_3)} & = \left(\pi, \dfrac{\pi}{2}, 0, 1, 0\right). \\
    \end{split}
\end{equation}
The samples highlighting variations in field orientation, directional shear, and magnitude, 
in the magnetic field, and the effects in the transport directions in the snapshots and the basis.

\begin{figure}
    \centering
    \includegraphics[width=0.32\linewidth]{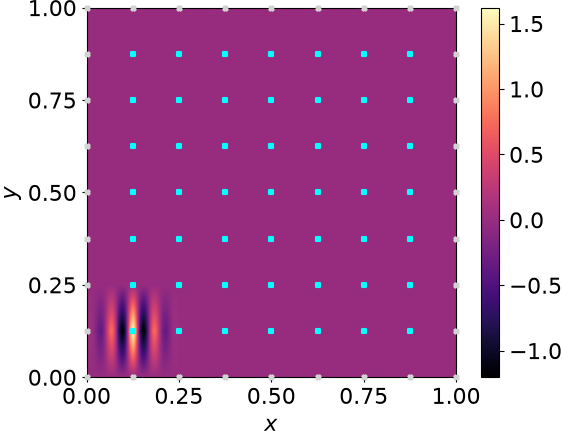}
    \includegraphics[width=0.32\linewidth]{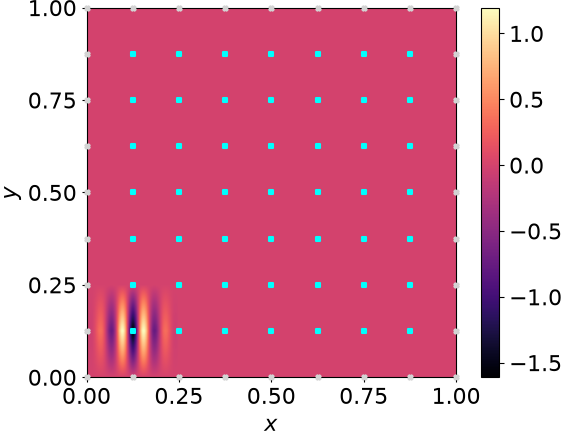}
    \includegraphics[width=0.32\linewidth]{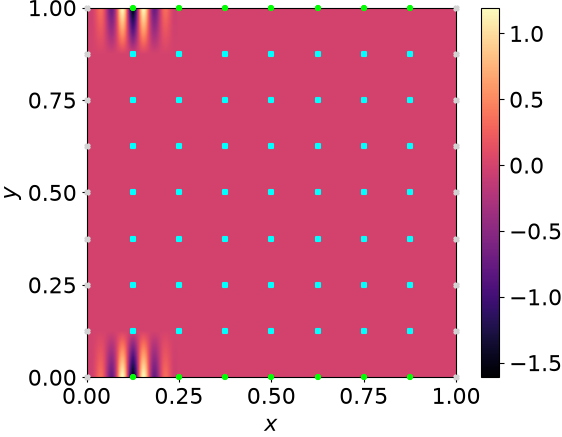}\\
    \includegraphics[width=0.32\linewidth]{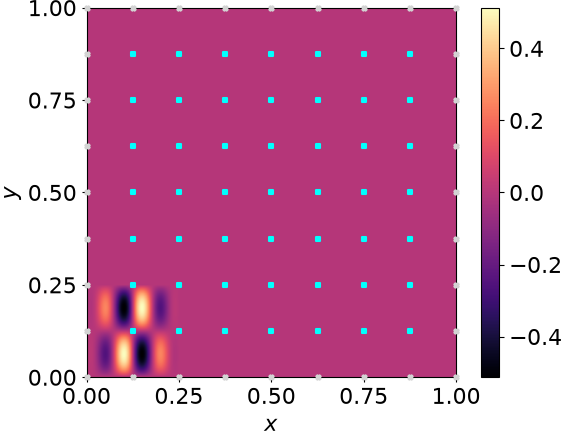}
    \includegraphics[width=0.32\linewidth]{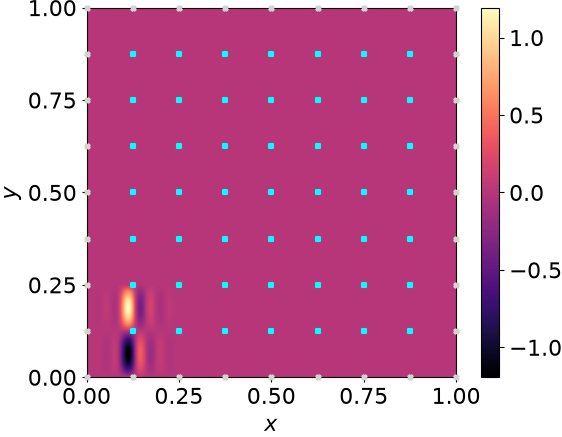}
    \includegraphics[width=0.32\linewidth]{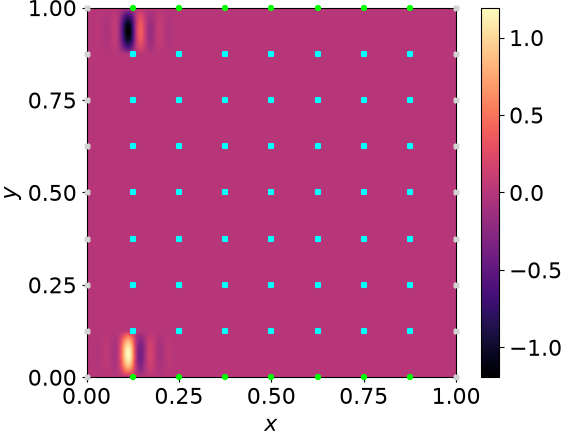}\\
    \includegraphics[width=0.32\linewidth]{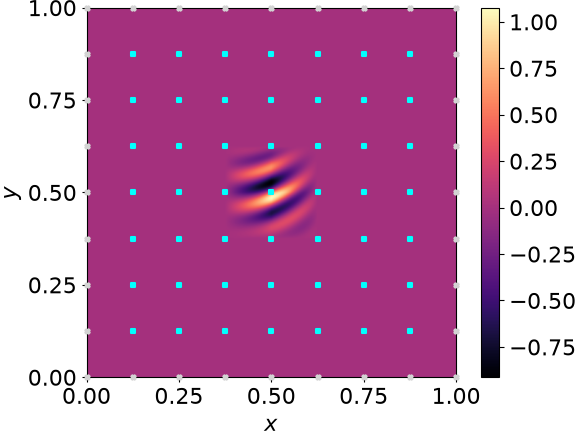}
    \includegraphics[width=0.32\linewidth]{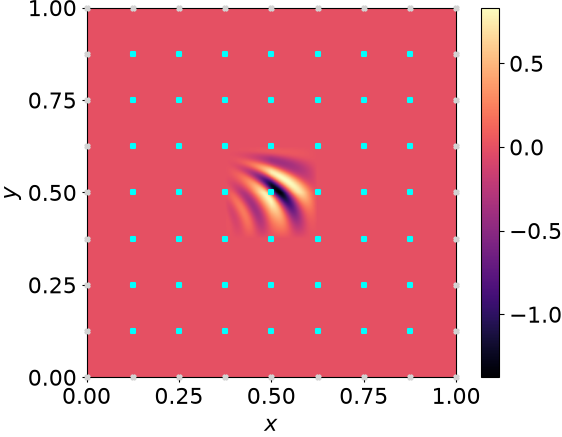}
    \includegraphics[width=0.32\linewidth]{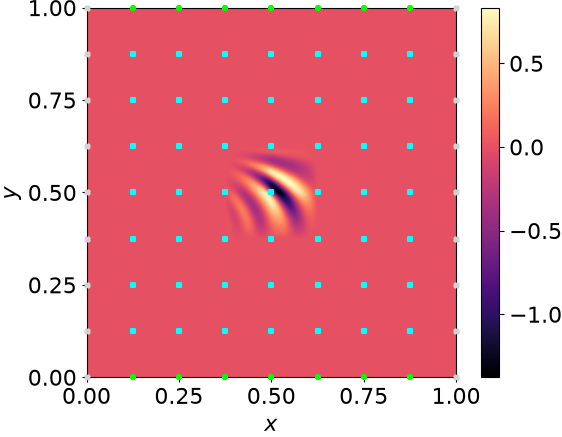}\\
    \includegraphics[width=0.32\linewidth]{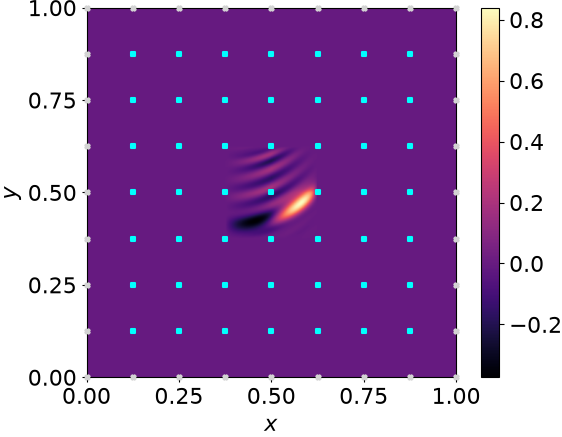}
    \includegraphics[width=0.32\linewidth]{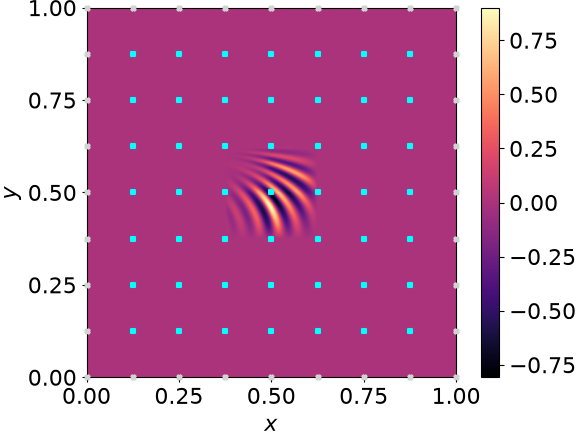}
    \includegraphics[width=0.32\linewidth]{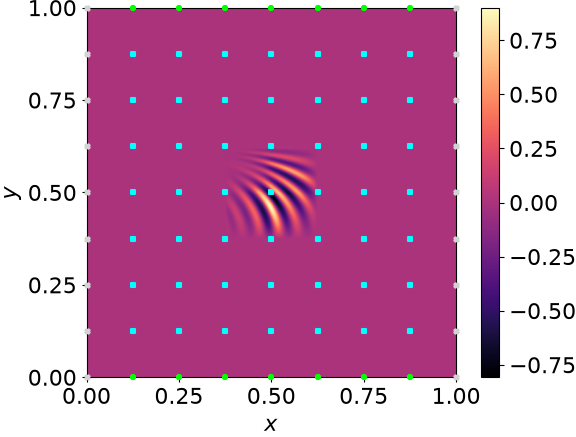}
    \caption{Representative samples of physical basis functions $\psi_{i_1,9}$ (1st row), $\psi_{i_1,17}$ (2nd row), $\psi_{i_2,9}$ (3rd row), and $\psi_{i_2,17}$ (4th row) for the single-null (left), double-null (middle), and magnetic island (right) problems.}
    \label{fig:assemble}
\end{figure}

In all benchmark test cases, the computation domain $\Omega = (0,1)^2$ is
partitioned into a coarse uniform mesh $\mathcal{T}_H$ of square elements with length $H$, with nodes denoted by $\{\mathbf{x}_i\}_{i\in\widetilde{\mathcal{I}}}$.
Let $\{\widetilde{\chi}_i\}_{i \in \widetilde{\mathcal{I}}} \subset H^1(\Omega)$ be the standard piecewise bilinear shape functions on $\mathcal{T}_H$,
satisfying the nodal property $\widetilde{\chi}_i(\mathbf{x}_j) = \delta_{ij}$.
We define the partition of unity $\{ \chi_i \}_{i \in \mathcal{I}}$, and the active index set $\mathcal{I}_{\text{DOF}}$
depending on the boundary condition prescribed.
\begin{itemize}
\item For the Dirichlet boundary condition in the single-null and double-null problems, 
\begin{itemize}
\item the partition of unity is indexed by the set of coarse nodes, i.e., $\mathcal{I} = \widetilde{\mathcal{I}}$, and consists of the standard bilinear shape functions on $\mathcal{T}_H$, i.e., 
\begin{equation}
\chi_i = \widetilde{\chi}_i \text{ for all } i \in \mathcal{I},
\label{eq:POU_Dirichlet}
\end{equation}
\item the active index set excludes the domain boundaries, i.e.,
\begin{equation}
\mathcal{I}_{\text{DOF}} = \{i \in \mathcal{I} : \mathbf{x}_i \notin \partial\Omega\}.
\label{eq:IDOF_Dirichlet}
\end{equation}
\end{itemize}
\item For the mixed Dirichlet-periodic boundary condition in the magnetic island problem, 
\begin{itemize}
\item the partition of unity is indexed by the set of coarse nodes excluding the top horizontal boundary, i.e., 
\begin{equation}
\mathcal{I} = \{i \in \widetilde{\mathcal{I}} : \mathbf{x}_i \notin [0,1] \times \{1\}\},
\end{equation}
and consists of modified shape functions for paired bottom boundary nodes and standard bilinear shape functions for all other nodes, i.e.,
\begin{equation}
\chi_i = 
\begin{cases}
\widetilde{\chi}_i & \text{ if } \mathbf{x}_i \notin [0,1] \times \{0, 1\}, \\
\widetilde{\chi}_i + \widetilde{\chi}_j & \text{ if } \mathbf{x}_i = (x, 0) \text{ and } \mathbf{x}_j = (x, 1), 
\end{cases}
\label{eq:POU_mixed}
\end{equation}
for all $i \in \mathcal{I}$. 
\item the active index set $\mathcal{I}_{\text{DOF}}$ excludes the vertical Dirichlet boundaries, i.e.,
\begin{equation}
\mathcal{I}_{\text{DOF}} = \{i \in \mathcal{I} : \mathbf{x}_i \notin \{0,1\} \times [0,1]\}.
\label{eq:IDOF_mixed}
\end{equation}
\end{itemize}
\end{itemize}
With the open cover $\mathcal{S}$ defined with $\omega_i = \text{supp}(\chi_i)^\circ$ for all $i \in \mathcal{I}$, 
these choices guarantee that $\{ \chi_i \}_{i \in \mathcal{I}}$ is a partition of unity subordinate to $\mathcal{S}$ and 
$\mathcal{I}_{\text{DOF}}$ counts all the degrees of freedom under the homogeneous essential boundary condition in the resolution of $\mathcal{T}_H$. 
In both cases, the mesh $\mathcal{T}_h$ is conforming to the cover $\{ \omega_i \}_{i \in \mathcal{I}}$, and 
each patch $\omega_i$ is isomorphic to the reference patch $\widehat{\omega} = (0,1)^2$. 
We remark that in the magnetic island problem, for patches centered along the bottom boundary nodes, i.e., $\mathbf{x}_i \in [0,1] \times \{ 0 \}$, 
the isomorphism $T_i$ reflects the underlying periodic magnetic topology and defines the local diffusion tensor $\widehat{\mathcal{K}}_i$ with the periodic information.
Figure~\ref{fig:assemble} demonstrates the physical basis functions $\psi_{i_1,9}$, $\psi_{i_1,17}$, $\psi_{i_2,9}$, and $\psi_{i_2,17}$ 
after identifying the library index $m_i^*$ in 2 representative patches $\omega_i$ with $i = i_1, i_2 \in \mathcal{I}_{\text{DOF}}$ for each of the benchmark problems. 
For the single-null and double-null problems, we select 
\begin{equation}
\mathbf{x}_{i_1} = \left(\dfrac{1}{8}, \dfrac{1}{8}\right), 
\mathbf{x}_{i_2} = \left(\dfrac{1}{8}, \dfrac{1}{2}\right).
\end{equation}
For the magnetic island problem, we select 
\begin{equation}
\mathbf{x}_{i_1} = \left(\dfrac{1}{8}, 0 \right), 
\mathbf{x}_{i_2} = \left(\dfrac{1}{8}, \dfrac{1}{2}\right).
\end{equation}
The cyan interior coarse nodes are always included in the active set $\mathcal{I}_{\text{DOF}}$, 
and the vertical boundary coarse nodes, painted in gray, are always excluded. 
The green coarse nodes on horizontal boundary in the last column represents the extra active indices 
in \eqref{eq:IDOF_mixed} specifically for the magnetic island problem, 
on which modified shape functions are defined in the partition of unity \eqref{eq:POU_mixed}.
The difference in the partition of unity and the active index set
between the magnetic island problem and the other two problems,
as stated in \eqref{eq:POU_Dirichlet}--\eqref{eq:IDOF_mixed}, is reflected 
in the physical basis functions supported in $\omega_{i_1}$ in the first two rows of Figure~\ref{fig:assemble}. 
On the other hand, the physical basis functions supported in $\omega_{i_2}$ in the last two rows of Figure~\ref{fig:assemble} reflect the similarity 
in the interior magnetic topology of the double-null configuration in Figure~\ref{fig:double_null_magnetic} and the magnetic island configuration in Figure~\ref{fig:magnetic_island_magnetic}.
While the number of local basis $r_i$ can be chosen independently for each patch, 
throughout all the numerical experiments in this work, we choose a constant number $r$ of local basis for all the patches. 

\subsection{PU-DDFEM solution and convergence}

In this study, the mesh parameter $H$ is chosen depending on the parallel diffusion coefficient $\kappa_{\parallel}$, 
namely $H = 2^{-(4+\log_{10} \sqrt{\kappa_{\parallel}})}$ for $\kappa_{\parallel} \in \{10^2, 10^4, 10^6\}$, 
to maintain the same coarsening ratio $H = 8h$ between $\mathcal{T}_h$ and $\mathcal{T}_H$. The prescribed frequency cutoff is $N_\text{freq} = 20$, and the number of local basis per patch is $r = 25$.

\begin{figure}[htbp!]
    \centering
    \includegraphics[width=0.32\linewidth]{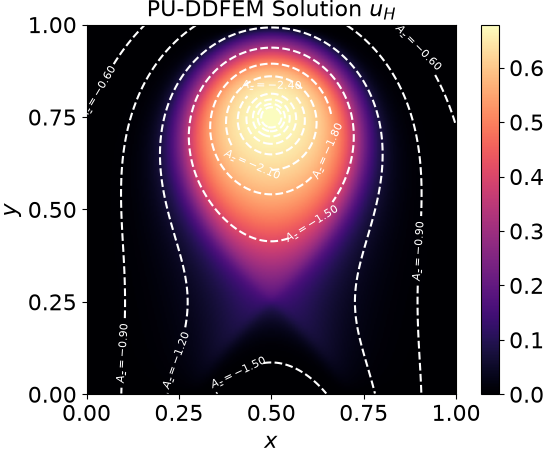}
    \includegraphics[width=0.32\linewidth]{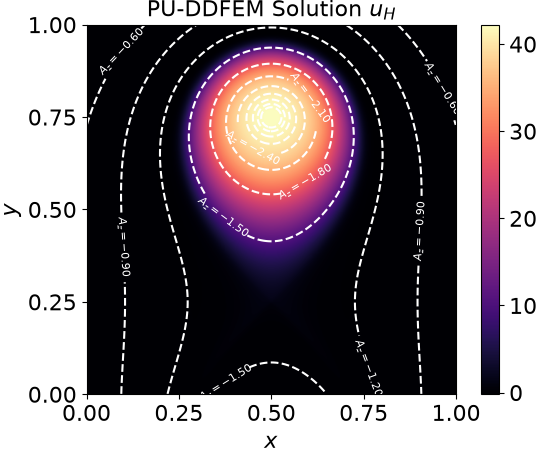}
    \includegraphics[width=0.32\linewidth]{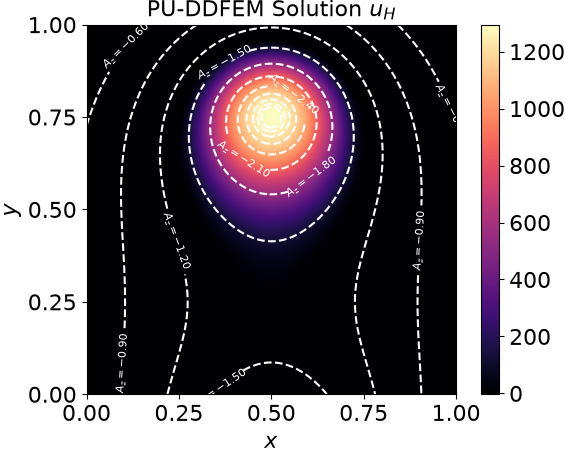}
    \caption{PU-DDFEM solution $u_H$ for the single-null problem with anisotropic ratio 
 $\kappa_{\parallel} = 10^2$ (left), $\kappa_{\parallel} = 10^4$ (middle), and $\kappa_{\parallel} = 10^6$ (right).}
    \label{fig:single_null_ddfem_solution}
\end{figure}

\begin{figure}[htbp!]
    \centering
    \includegraphics[width=0.32\linewidth]{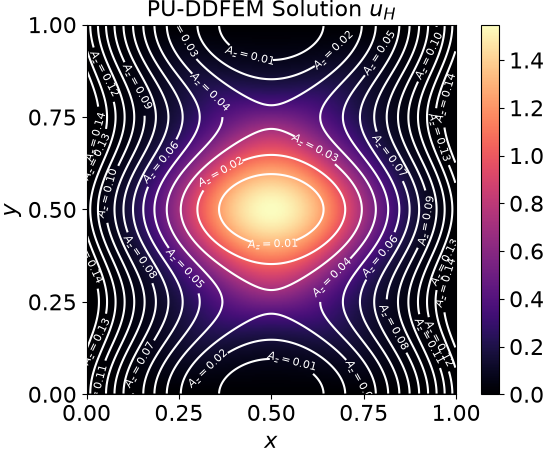}
    \includegraphics[width=0.32\linewidth]{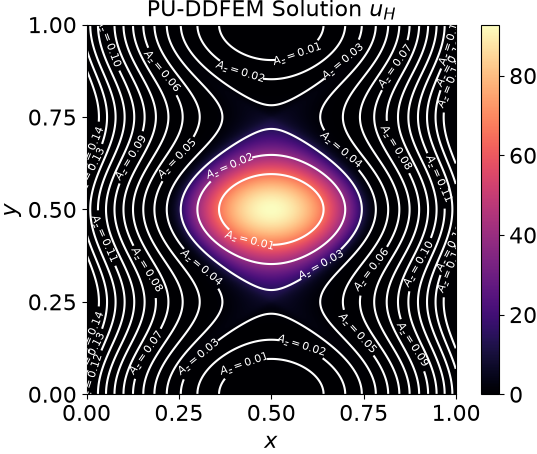}
    \includegraphics[width=0.32\linewidth]{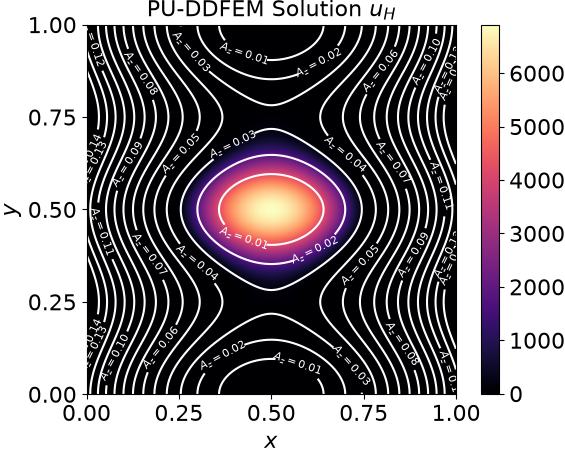}
    \caption{PU-DDFEM solution $u_H$ for the double-null problem with anisotropic ratio $\kappa_{\parallel} = 10^2$ (left), $\kappa_{\parallel} = 10^4$ (middle), and $\kappa_{\parallel} = 10^6$ (right).}
    \label{fig:double_null_ddfem_solution}
\end{figure}

\begin{figure}[htbp!]
    \centering
    \includegraphics[width=0.32\linewidth]{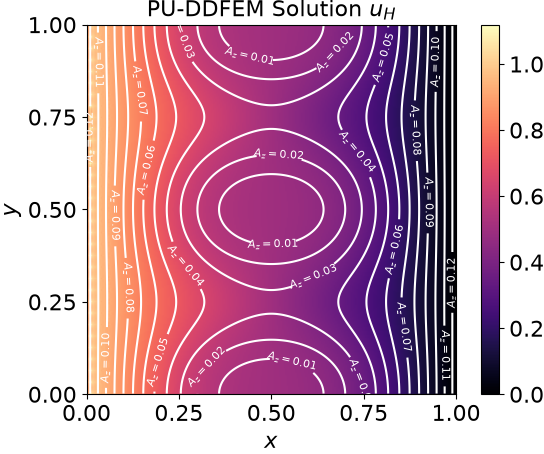}
    \includegraphics[width=0.32\linewidth]{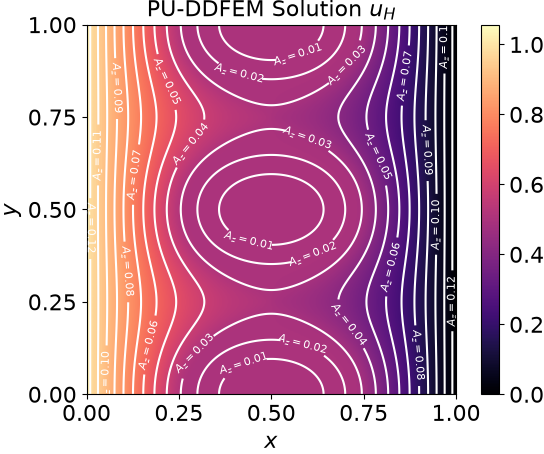}
    \includegraphics[width=0.32\linewidth]{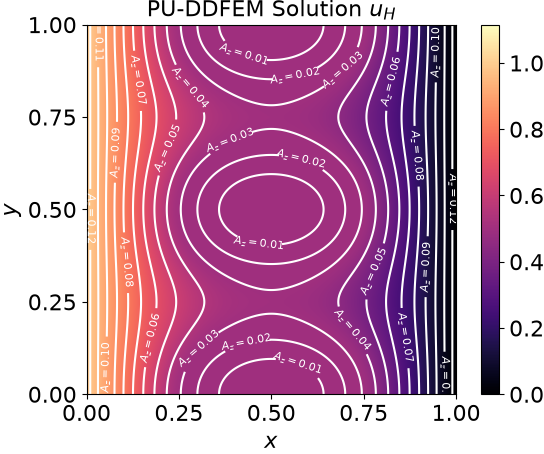}
    \caption{PU-DDFEM solution $u_H$ for the magnetic island problem with anisotropic ratio $\kappa_{\parallel} = 10^2$ (left), $\kappa_{\parallel} = 10^4$ (middle), and $\kappa_{\parallel} = 10^6$ (right).}
    \label{fig:magnetic_island_ddfem_solution}
\end{figure}

Figures~\ref{fig:single_null_ddfem_solution}, \ref{fig:double_null_ddfem_solution}, and \ref{fig:magnetic_island_ddfem_solution} 
shows the PU-DDFEM solution in \eqref{eq:PU-DDFEM-sol} 
for the three benchmark problems with anisotropic ratio $\kappa_{\parallel} \in \{10^2, 10^4, 10^6\}$. 
It can be seen that the solutions visually agree with the reference FEM solution in 
Figures~\ref{fig:single_null_fem_solution}, \ref{fig:double_null_fem_solution}, and \ref{fig:magnetic_island_fem_solution} 
well. 

Table~\ref{tab:ddfem_r25_results} compares the metrics for the three benchmark problems with anisotropic ratio $\kappa_{\parallel} \in \{10^2, 10^4, 10^6\}$. Overall, the relative $L^2$ errors remain within the order of $10^{-5}$ to $10^{-2}$ across all benchmark problems, indicating that PU-DDFEM maintains reasonable accuracy over a wide range of anisotropic ratios. Meanwhile, the speed-up $\tau^{-1}$ generally increases to $> 3000 \times$ with $\kappa_{\parallel}$, suggesting that the method becomes increasingly advantageous as the anisotropy strengthens. This trend indicates the potential of PU-DDFEM for problems with even higher anisotropic ratios. The performance also depends on the underlying magnetic geometry, with the magnetic island problem exhibiting a more pronounced increase in computational cost at high anisotropy.

\begin{table}[htbp]
    \centering
    \begin{tabular}{lc|cc}
        \toprule
        Problem & $\kappa_{\parallel}$ & $\varepsilon_{L^2}$ & $\tau^{-1}$ \\
        \midrule
        \multirow{3}{*}{Single-null} 
            & $10^2$ & $1.41 \times 10^{-4}$ & 656.61 \\
            & $10^4$ & $6.44 \times 10^{-3}$ & 3441.87 \\
            & $10^6$ & $1.74 \times 10^{-2}$ & 6387.08 \\
        \midrule
        \multirow{3}{*}{Double-null} 
            & $10^2$ & $4.72 \times 10^{-5}$ & 228.40 \\
            & $10^4$ & $1.15 \times 10^{-3}$ & 1109.58 \\
            & $10^6$ & $6.71 \times 10^{-3}$ & 3563.27 \\
        \midrule
        \multirow{3}{*}{Magnetic island} 
            & $10^2$ & $9.63 \times 10^{-3}$ & 158.52 \\
            & $10^4$ & $5.19 \times 10^{-3}$ & 1320.67 \\
            & $10^6$ & $8.33 \times 10^{-3}$ & 12369.59 \\
        \bottomrule
    \end{tabular}
    \caption{Relative $L^2$ error and speed-up of PU-DDFEM for the benchmark problems}
    \label{tab:ddfem_r25_results}
\end{table}

To evaluate the impact of the local basis size $r$ on numerical accuracy and computational efficiency, we examine the convergence behavior across three distinct magnetic configurations, single-null, double-null, and magnetic island, under varying anisotropic ratios $\kappa_{\parallel} \in \{10^2, 10^4, 10^6\}$.

\begin{figure}[htbp!]
    \centering
    \includegraphics[width=0.32\linewidth]{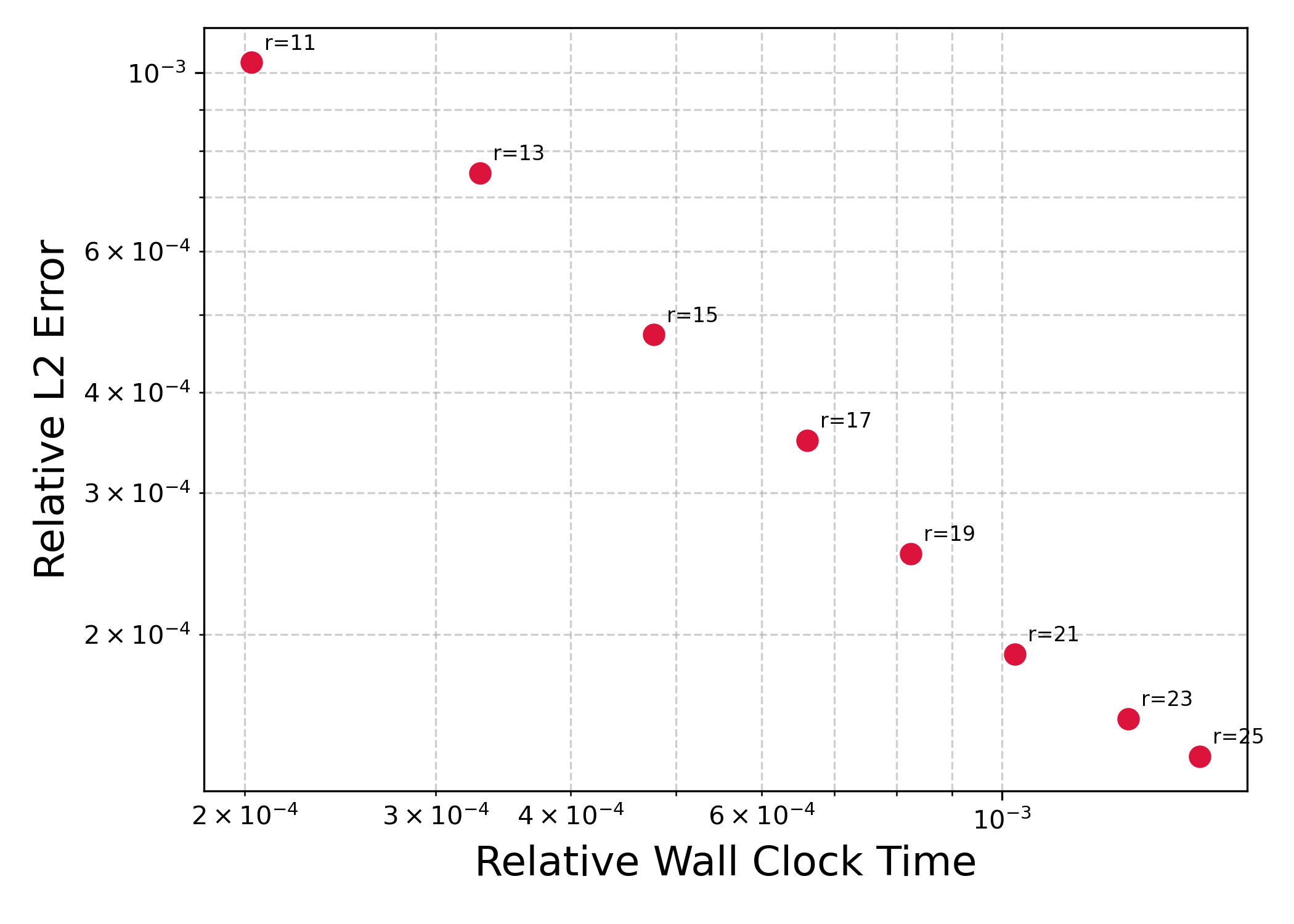}
    \includegraphics[width=0.32\linewidth]{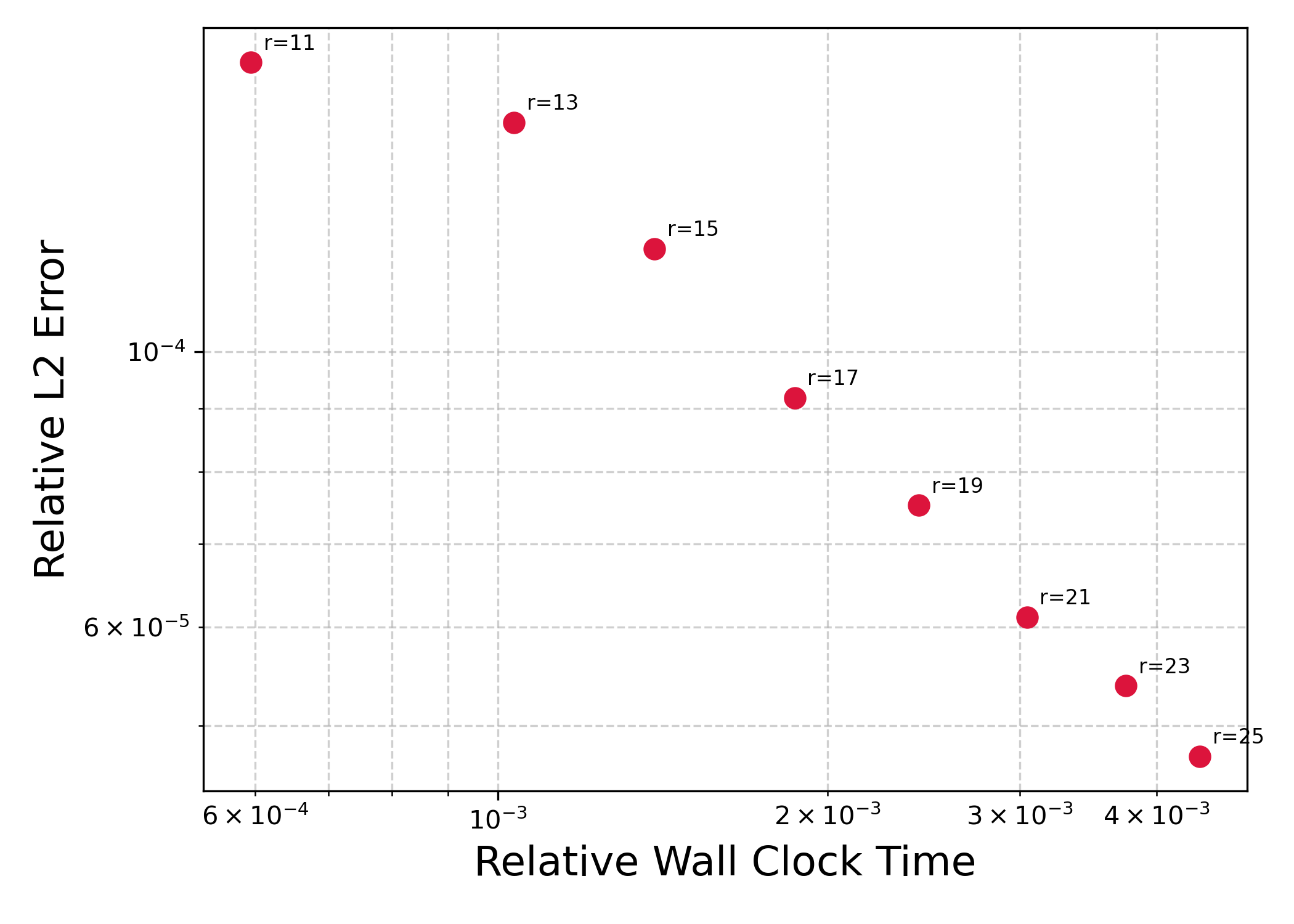}
    \includegraphics[width=0.32\linewidth]{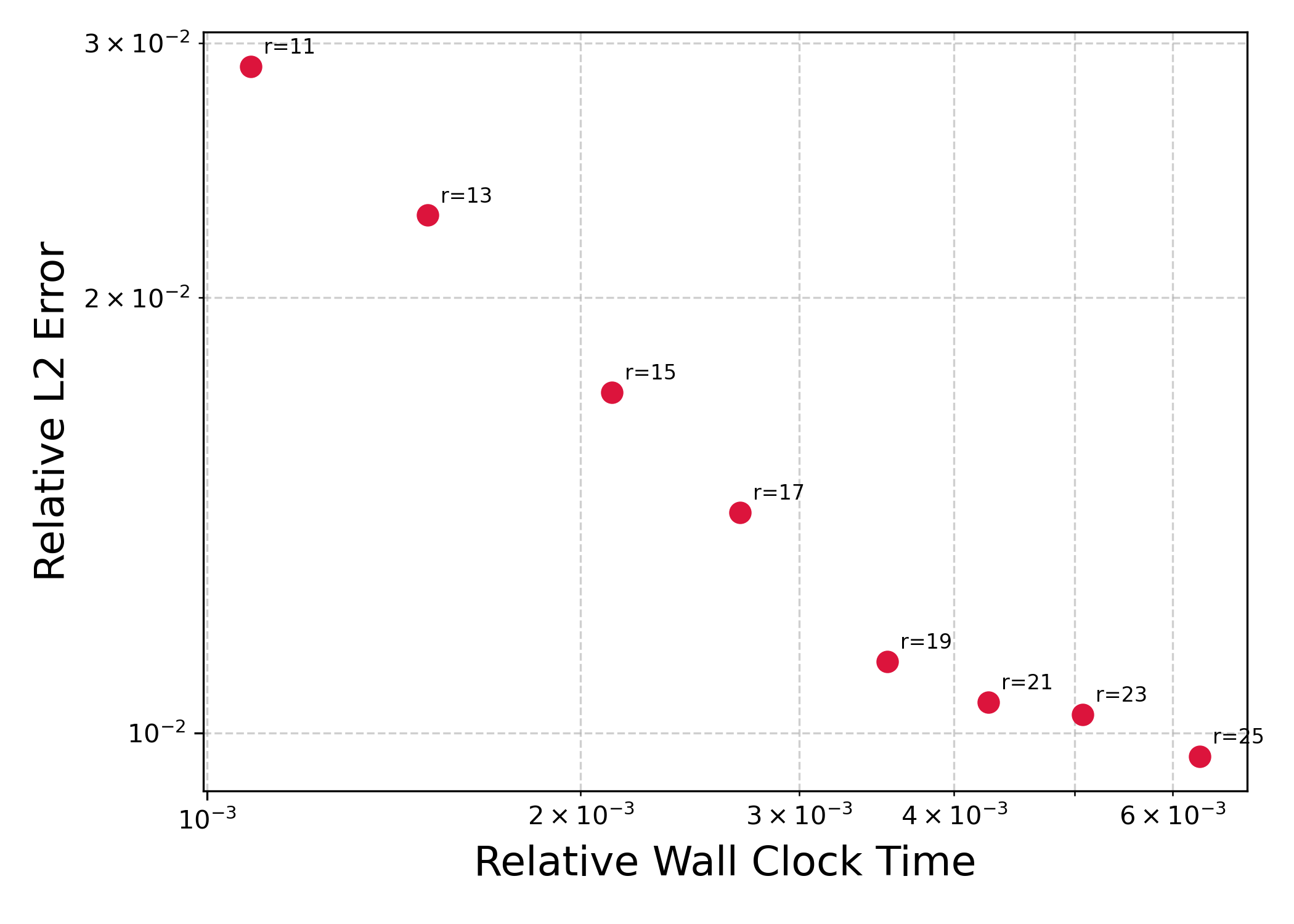}
    \caption{Pareto optimal of competing interests of $L^2$ error and relative wall clock time for the single-null problem (left), double-null problem (middle), and magnetic island problem (right) with anisotropic ratio $\kappa_{\parallel} = 10^2$.}
    \label{fig:anis2_performance}
\end{figure}

For a mild anisotropy ratio of $\kappa_{\parallel} = 10^2$, all three configurations exhibit rapid error decay with modest basis sizes, as illustrated in Fig.~\ref{fig:anis2_performance}. 
The results demonstrate a consistent reduction in the relative $L^2$ error as the local basis size $r$ increases across all three benchmark configurations. The single-null and double-null cases exhibit particularly high accuracy, with errors reaching the orders of $10^{-4}$ and $10^{-5}$, respectively, for the largest basis sizes considered. In contrast, the magnetic island benchmark is more challenging, with the relative error remaining on the order of $10^{-2}$ despite increasing the local basis size. This indicates that the approximation quality depends strongly on the underlying magnetic geometry, with the magnetic island configuration exhibiting a more pronounced accuracy limitation. The computational speed-up decreases monotonically as the local basis size increases, reflecting the additional computational cost associated with enriching the reduced approximation space. Nevertheless, substantial acceleration is maintained throughout the tested range. The speed-up factors remain on the order of $10^2$--$10^3$ for the double-null and magnetic island configurations and reach the order of $10^3$--$10^4$ for the single-null case at smaller basis sizes. Overall, these results illustrate that increasing the local basis size provides a systematic improvement in accuracy while retaining significant computational savings relative to the reference discretization.

\begin{figure}[htbp!]
    \centering
    \includegraphics[width=0.32\linewidth]{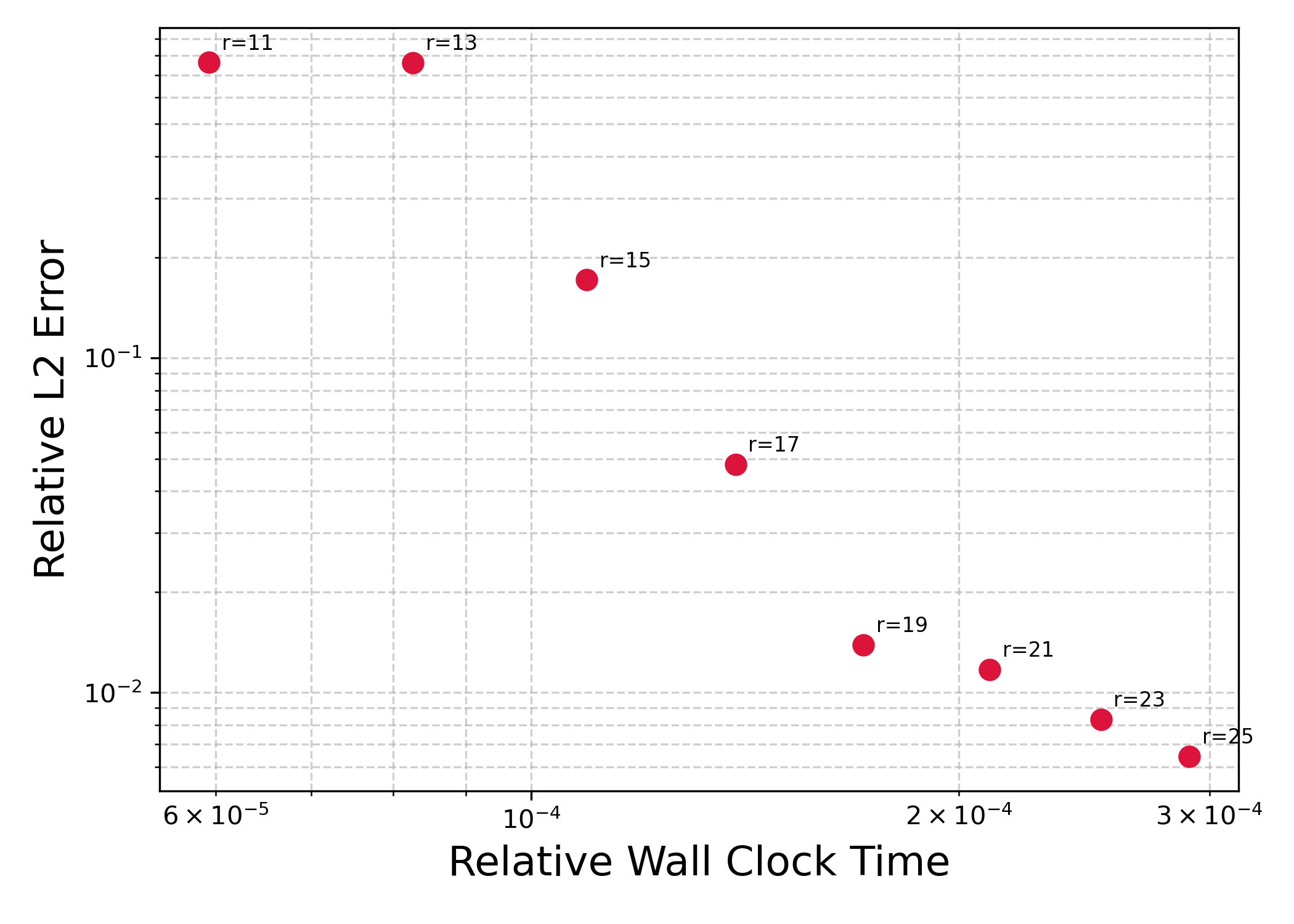}
    \includegraphics[width=0.32\linewidth]{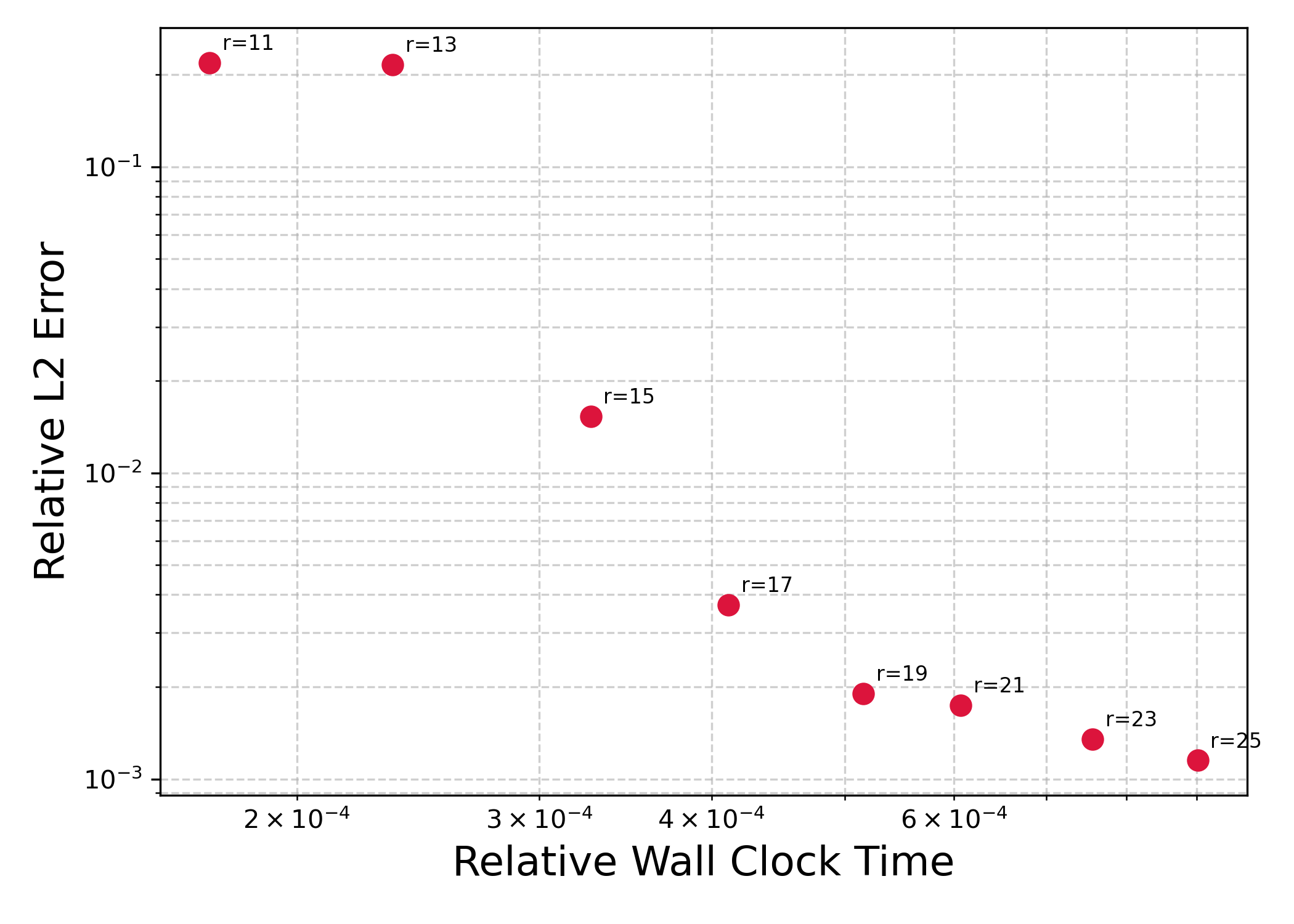}
    \includegraphics[width=0.32\linewidth]{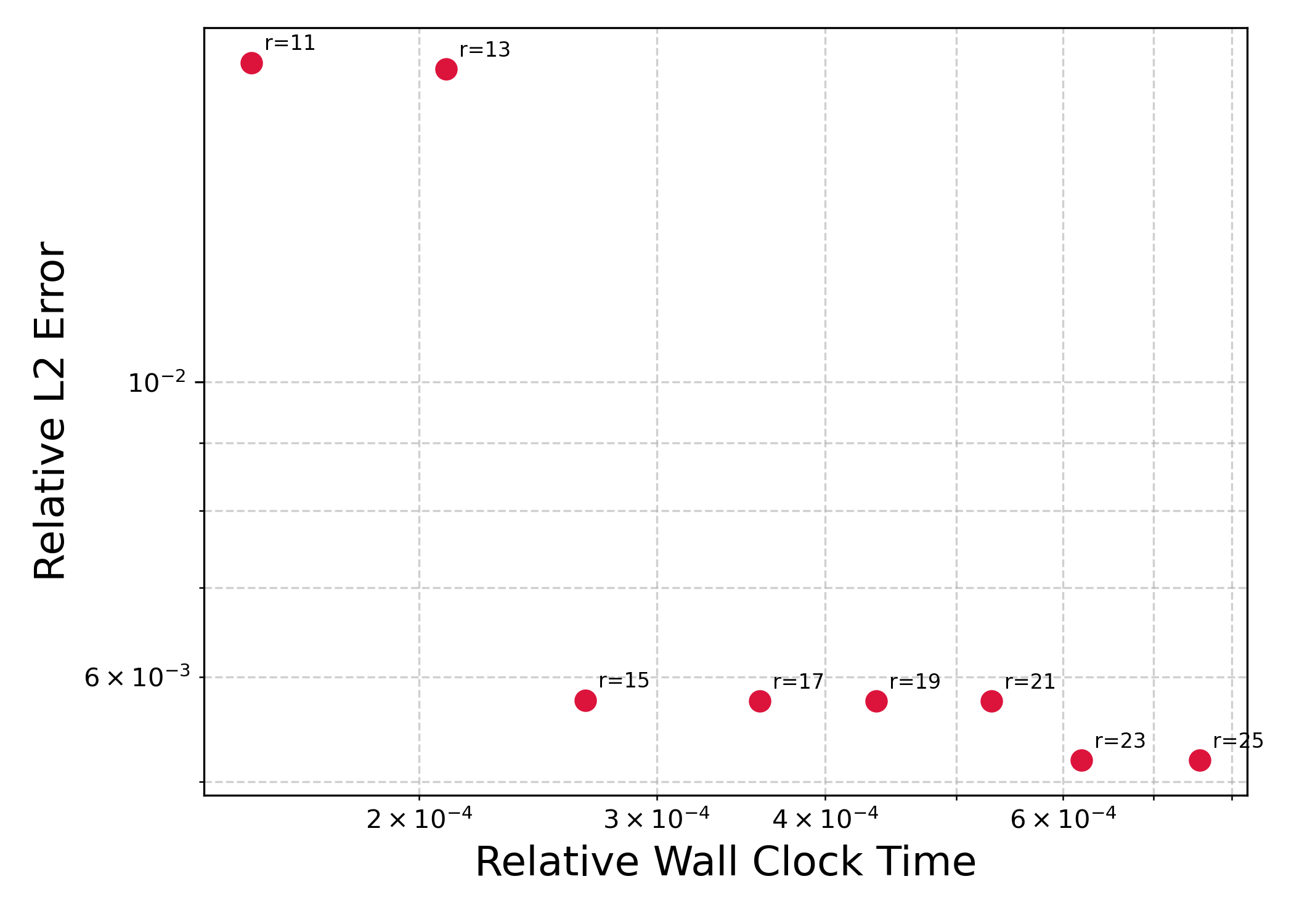}
    \caption{Pareto optimal of competing interests of $L^2$ error and relative wall clock time for the single-null problem (left), double-null problem (middle), and magnetic island problem (right) with anisotropic ratio $\kappa_{\parallel} = 10^4$.}
    \label{fig:anis4_performance}
\end{figure}

As the anisotropy increases to $\kappa_{\parallel} = 10^4$, as shown in Figure~\ref{fig:anis4_performance}, the approximation becomes more sensitive to the local basis size. In both the single-null and double-null benchmarks, the error remains relatively large for the smallest basis sizes, followed by a substantial reduction as the basis is enriched. The resulting errors reach the order of $10^{-3}$ for both configurations at the largest basis sizes. This behavior suggests that a sufficiently rich local approximation space is essential for accurately resolving the solution under strong anisotropy.
The magnetic island benchmark displays a different convergence pattern. Although the relative error decreases from the order of $10^{-2}$ to several times $10^{-3}$ as the local basis size increases, the improvement becomes limited beyond moderate basis enrichment. This plateau indicates that increasing the number of local basis functions alone may not fully eliminate the approximation error for this geometry. Despite these differences in accuracy, all three configurations exhibit substantial computational speed-ups, ranging from several hundred to several thousand, depending on the benchmark and basis size. As expected, the speed-up decreases with increasing $r$, highlighting the trade-off between reduced-order accuracy and computational efficiency.

\begin{figure}[htbp!]
    \centering
    \includegraphics[width=0.32\linewidth]{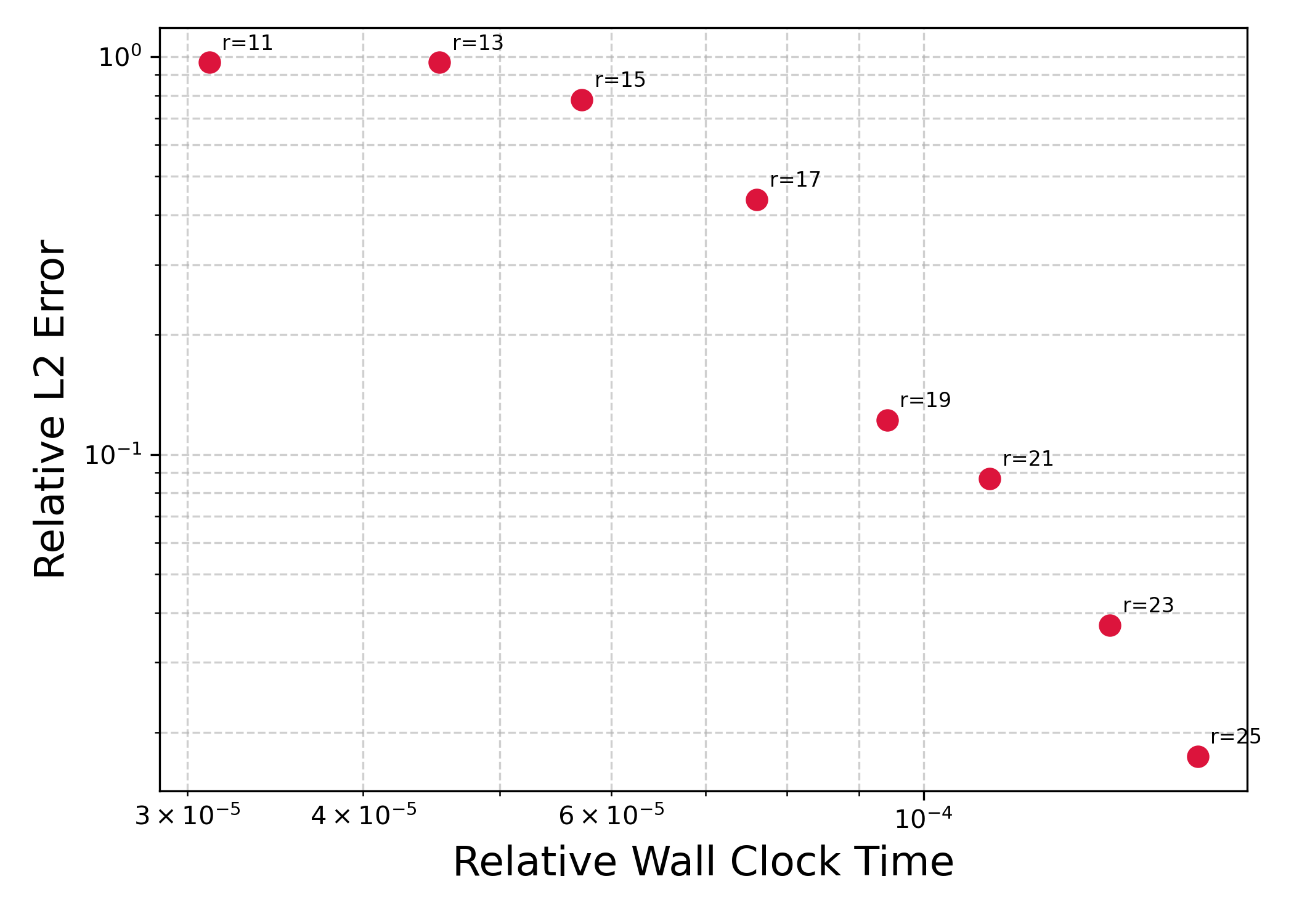}
    \includegraphics[width=0.32\linewidth]{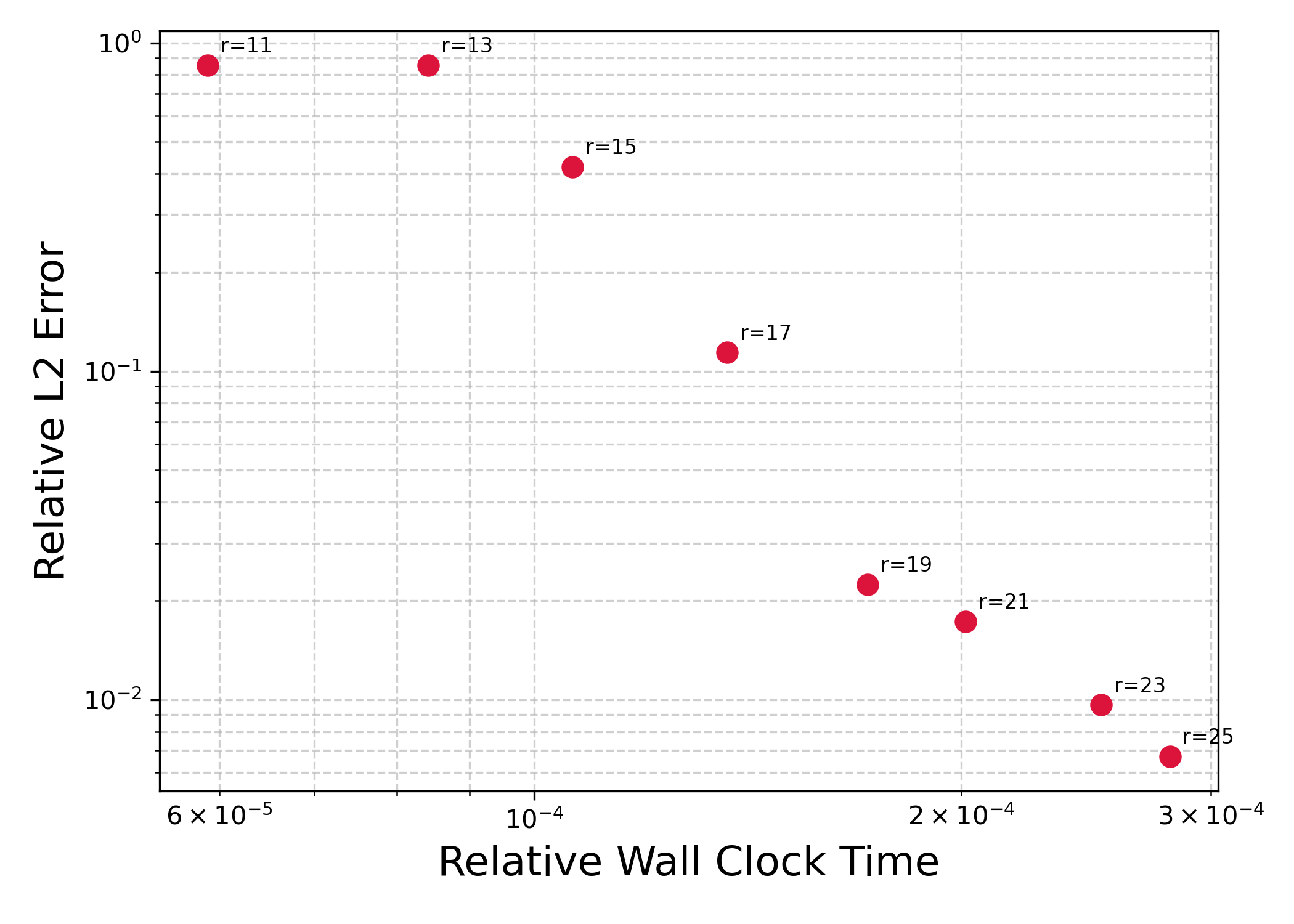}
    \includegraphics[width=0.32\linewidth]{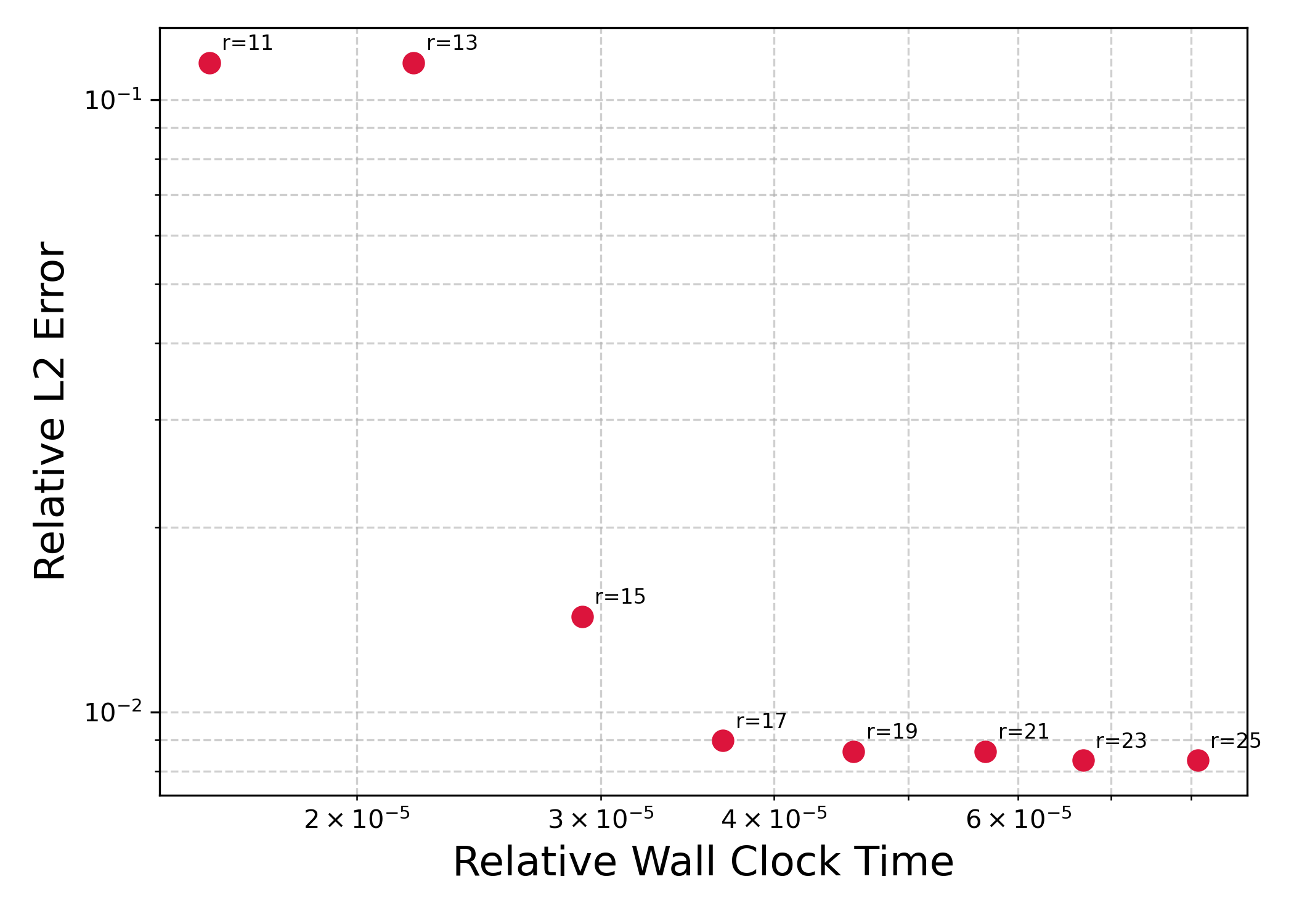}
    \caption{Pareto optimal of competing interests of $L^2$ error and relative wall clock time for the single-null problem (left), double-null problem (middle), and magnetic island problem (right) with anisotropic ratio $\kappa_{\parallel} = 10^6$.}
    \label{fig:anis6_performance}
\end{figure}

Under high anisotropy $\kappa_{\parallel} = 10^6$, the necessity of a sufficient local basis size becomes even more obvious, as shown in Figure~\ref{fig:anis6_performance}. 
For both the single-null and double-null configurations, the relative $L^2$ error remains close to order unity for the smallest basis sizes, indicating that the corresponding reduced spaces are insufficient to capture the solution accurately. As the local basis is enriched, the error decreases sharply by several orders of magnitude, eventually reaching the order of $10^{-2}$ for the single-null case and $10^{-3}$--$10^{-2}$ for the double-null case. These results emphasize the importance of local basis enrichment when the anisotropy is particularly strong.
The magnetic island benchmark exhibits comparatively favorable accuracy at smaller basis sizes, with errors initially on the order of $10^{-1}$. A substantial improvement occurs upon increasing the basis size, after which the error stabilizes at several times $10^{-3}$. The observed plateau suggests that further enrichment yields diminishing returns and that the remaining error may be influenced by factors beyond the local basis dimension. Despite the severe anisotropy and the varying convergence behavior, the method maintains exceptionally large computational speed-ups, ranging from several thousand to tens of thousands. The results therefore demonstrate that the proposed approach retains considerable computational advantages even in highly anisotropic settings, although achieving acceptable accuracy requires a geometry- and parameter-dependent choice of local basis size.

\section{Conclusion}
\label{sec:conclusion}

In this work, we presented PU-DDFEM as a data-driven bottom-up local model reduction framework for solving anisotropic diffusion problems in strongly magnetized fusion plasmas. By parametrizing local magnetic field topologies, including orientation, shear, and curvature, over canonical reference patches, we constructed a reusable offline reduced basis library using complete trigonometric boundary excitations and localized generalized spectral decompositions. In the fast online phase, local subdomains select optimal pre-trained spectral modes via $L^2$-tensor matching, which are seamlessly glued into a globally $H^1$-conforming multiscale approximation space using a partition of unity.

Numerical benchmark evaluations across three representative plasma topologies, single-null, double-null, and magnetic island configurations, demonstrate the efficiency and accuracy of our approach under strong anisotropy. The method successfully resolves complex internal boundary layers, cross-wind flux behavior, and field line reconnection dynamics. The proposed PU-DDFEM scheme achieves $L^2$ errors under $2\%$ while delivering substantial computational speed-ups of $> 3000 \times$ over high-fidelity full-order finite element solves. 

Several promising directions remain for future investigation. First, a systematic comparison with alternative approaches for local reduced-basis construction and coupling \cite{efendiev2013generalized-dg,chung2015mixed,maalqvist2020numerical,altmann2021numerical,chung2018constraint,cheung2020constraint,cheung2022iterative} would help quantify the advantages and limitations of PU-DDFEM in terms of accuracy, computational cost, and robustness. Second, incorporating online adaptivity \cite{chung2014adaptive,chung2017online,pun2023online} into the local basis selection and enrichment process could further improve efficiency by dynamically allocating reduced degrees of freedom according to local approximation errors, magnetic field variations, and evolving boundary-layer structures. Finally, extending PU-DDFEM with magnetohydrodynamic (MHD) models would provide a pathway toward self-consistent simulations of plasma transport and magnetic field evolution, enabling the treatment of more realistic fusion scenarios in which the magnetic geometry and anisotropic transport are mutually coupled. These extensions would advance PU-DDFEM toward an adaptive and scalable computational framework for high-fidelity plasma simulation.

\section*{Acknowledgments}\label{sec:disc}
This work was supported in part by the U.S. Department of Energy,
Office of Science, Office of Advanced Scientific Computing Research,
Scientific Discovery through Advanced Computing (SciDAC) program through
the LEADS SciDAC Institute under Project Number SCW1933 at
Lawrence Livermore National Laboratory (LLNL). 
LLNL is operated by Lawrence Livermore National Security LLC, for the
U.S. Department of Energy, National Nuclear Security Administration
under Contract DE-AC52-07NA27344.
IM release number: LLNL-JRNL-2025865.

\section*{Disclaimer}
This document was prepared as an account of work sponsored by an
agency of the United States government.
Neither the United States government nor Lawrence Livermore
National Security LLC nor any of their employees makes any warranty,
expressed or implied, or assumes any legal liability or responsibility
for the accuracy, completeness, or usefulness of any information,
apparatus, product, or process disclosed, or represents that its use
would not infringe privately owned rights.
Reference herein to any specific commercial product, process, or
service by trade name, trademark, manufacturer, or otherwise does not
necessarily constitute or imply its endorsement, recommendation,
or favoring by the United States government or
Lawrence Livermore National Security LLC.
The views and opinions of authors expressed herein do not necessarily
state or reflect those of the United States government or
Lawrence Livermore National Security LLC
and shall not be used for advertising or product endorsement purposes.

\bibliographystyle{plain}
\bibliography{references}

@article{spitzer1953transport,
  title={Transport phenomena in a completely ionized gas},
  author={Spitzer Jr, Lyman and H{\"a}rm, Richard},
  journal={Physical Review},
  volume={89},
  number={5},
  pages={977},
  year={1953},
  publisher={APS}
}

@article{braginskii1965transport,
  title={Transport processes in a plasma},
  author={Braginskii, SI},
  journal={Reviews of plasma physics},
  volume={1},
  pages={205},
  year={1965}
}

@article{gunter2005modelling,
  title={Modelling of heat transport in magnetised plasmas using non-aligned coordinates},
  author={G{\"u}nter, S and Yu, Q and Kr{\"u}ger, J and Lackner, K},
  journal={Journal of Computational Physics},
  volume={209},
  number={1},
  pages={354--370},
  year={2005},
  publisher={Elsevier}
}

@article{vogl2023mesh,
  title={Mesh refinement for anisotropic diffusion in magnetized plasmas},
  author={Vogl, Christopher J and Joseph, Ilon and Holec, Milan},
  journal={Computers \& Mathematics with Applications},
  volume={145},
  pages={159--174},
  year={2023},
  publisher={Elsevier}
}

@article{melenk1996partition,
  title={The partition of unity finite element method: basic theory and applications},
  author={Melenk, Jens M and Babuska, Ivo},
  journal={Computer methods in applied mechanics and engineering},
  volume={139},
  number={1-4},
  pages={289--314},
  year={1996},
  publisher={Elsevier}
}

@article{babuska1997partition,
  title={The partition of unity method},
  author={Babuska, Ivo and Melenk, Jens M},
  journal={International journal for numerical methods in engineering},
  volume={40},
  number={4},
  pages={727--758},
  year={1997},
  publisher={Wiley Online Library}
}

@article{babuvska2004generalized,
  title={Generalized finite element methods—main ideas, results and perspective},
  author={Babuska, Ivo and Banerjee, Uday and Osborn, John E},
  journal={International Journal of Computational Methods},
  volume={1},
  number={01},
  pages={67--103},
  year={2004},
  publisher={World Scientific}
}

@book{hesthaven2016certified,
  title={Certified reduced basis methods for parametrized partial differential equations},
  author={Hesthaven, Jan S and Rozza, Gianluigi and Stamm, Benjamin},
  volume={590},
  year={2016},
  publisher={Springer}
}

@book{quarteroni2015reduced,
  title={Reduced Basis Methods for Partial Differential Equations: An Introduction},
  author={Quarteroni, Alfio and Manzoni, Andrea and Negri, Federico},
  volume={92},
  year={2015},
  publisher={Springer}
}

@article{huynh2013static,
  title={A static condensation reduced basis element method: approximation and a posteriori error estimation},
  author={Huynh, Dinh Bao Phuong and Knezevic, David J and Patera, Anthony T},
  journal={ESAIM: Mathematical Modelling and Numerical Analysis},
  volume={47},
  number={1},
  pages={213--251},
  year={2013},
  publisher={EDP Sciences}
}

@article{hoang2021domain,
  title={Domain-decomposition least-squares {Petrov--Galerkin (DD-LSPG)} nonlinear model reduction},
  author={Hoang, Chi and Choi, Youngsoo and Carlberg, Kevin},
  journal={Computer methods in applied mechanics and engineering},
  volume={384},
  pages={113997},
  year={2021},
  publisher={Elsevier}
}

@article{smetana2016optimal,
  title={Optimal local approximation spaces for component-based static condensation procedures},
  author={Smetana, Kathrin and Patera, Anthony T},
  journal={SIAM Journal on Scientific Computing},
  volume={38},
  number={5},
  pages={A3318--A3356},
  year={2016},
  publisher={SIAM}
}

@article{ebrahimi2024hyperreduced,
  title={A hyperreduced reduced basis element method for reduced-order modeling of component-based nonlinear systems},
  author={Ebrahimi, Mehran and Yano, Masayuki},
  journal={Computer Methods in Applied Mechanics and Engineering},
  volume={431},
  pages={117254},
  year={2024},
  publisher={Elsevier}
}

@article{zanardi2024scalable,
  title={Scalable nonlinear manifold reduced order model for dynamical systems},
  author={Zanardi, Ivan and Diaz, Alejandro N and Chung, Seung Whan and Panesi, Marco and Choi, Youngsoo},
  journal={Machine Learning and the Physical Sciences Workshop, NeurIPS 2024},
  year={2024}
}

@article{chung2024scalable,
  title={Scalable physics-guided data-driven component model reduction for steady Navier-Stokes flow},
  author={Chung, Seung Whan and Choi, Youngsoo and Roy, Pratanu and Roy, Thomas and Lin, Tiras Y and Nguyen, Du T and Hahn, Christopher and Duoss, Eric B and Baker, Sarah E},
  journal={Machine Learning and the Physical Sciences Workshop, NeurIPS 2024},
  year={2024}
}

@article{chung2024scaled,
  title={Scaled-up prediction of steady Navier-Stokes equation with component reduced order modeling},
  author={Chung, Seung Whan and Choi, Youngsoo and Roy, Pratanu and Roy, Thomas and Lin, Tiras Y and Nguyen, Du T and Hahn, Christopher and Duoss, Eric B and Baker, Sarah E},
  journal={arXiv preprint arXiv:2410.21534},
  year={2024}
}

@article{eftang2013port,
  title={Port reduction in parametrized component static condensation: approximation and a posteriori error estimation},
  author={Eftang, Jens L and Patera, Anthony T},
  journal={International Journal for Numerical Methods in Engineering},
  volume={96},
  number={5},
  pages={269--302},
  year={2013},
  publisher={Wiley Online Library}
}

@article{choi2019space,
  title={Space--time least-squares {Petrov--Galerkin} projection for nonlinear model reduction},
  author={Choi, Youngsoo and Carlberg, Kevin},
  journal={SIAM Journal on Scientific Computing},
  volume={41},
  number={1},
  pages={A26--A58},
  year={2019},
  publisher={SIAM}
}

@article{chung2024train,
  title={Train small, model big: Scalable physics simulators via reduced order modeling and domain decomposition},
  author={Chung, Seung Whan and Choi, Youngsoo and Roy, Pratanu and Moore, Thomas and Roy, Thomas and Lin, Tiras Y and Nguyen, Du T and Hahn, Christopher and Duoss, Eric B and Baker, Sarah E},
  journal={Computer Methods in Applied Mechanics and Engineering},
  volume={427},
  pages={117041},
  year={2024},
  publisher={Elsevier}
}

@article{carlberg2018conservative,
  title={Conservative model reduction for finite-volume models},
  author={Carlberg, Kevin and Choi, Youngsoo and Sargsyan, Syuzanna},
  journal={Journal of Computational Physics},
  volume={371},
  pages={280--314},
  year={2018},
  publisher={Elsevier}
}

@article{lee2020model,
  title={Model reduction of dynamical systems on nonlinear manifolds using deep convolutional autoencoders},
  author={Lee, Kookjin and Carlberg, Kevin T},
  journal={Journal of Computational Physics},
  volume={404},
  pages={108973},
  year={2020},
  publisher={Elsevier}
}

@article{wentland2025role,
  title={The role of interface boundary conditions and sampling strategies for Schwarz-based coupling of projection-based reduced order models},
  author={Wentland, Christopher R and Rizzi, Francesco and Barnett, Joshua L and Tezaur, Irina K},
  journal={Journal of Computational and Applied Mathematics},
  volume={465},
  pages={116584},
  year={2025},
  publisher={Elsevier}
}

@article{diaz2024fast,
  title={A fast and accurate domain decomposition nonlinear manifold reduced order model},
  author={Diaz, Alejandro N and Choi, Youngsoo and Heinkenschloss, Matthias},
  journal={Computer Methods in Applied Mechanics and Engineering},
  volume={425},
  pages={116943},
  year={2024},
  publisher={Elsevier}
}

@article{eftang2014port,
  title={A port-reduced static condensation reduced basis element method for large component-synthesized structures: approximation and a posteriori error estimation},
  author={Eftang, Jens L and Patera, Anthony T},
  journal={Advanced Modeling and Simulation in Engineering Sciences},
  volume={1},
  pages={1--49},
  year={2014},
  publisher={Springer}
}

@article{huynh2013static2,
  title={A static condensation reduced basis element method: Complex problems},
  author={Huynh, Dinh Bao Phuong and Knezevic, David J and Patera, Anthony T},
  journal={Computer Methods in Applied Mechanics and Engineering},
  volume={259},
  pages={197--216},
  year={2013},
  publisher={Elsevier}
}

@article{mcbane2022stress,
  title={Stress-constrained topology optimization of lattice-like structures using component-wise reduced order models},
  author={McBane, Sean and Choi, Youngsoo and Willcox, Karen},
  journal={Computer Methods in Applied Mechanics and Engineering},
  volume={400},
  pages={115525},
  year={2022},
  publisher={Elsevier}
}

@article{mcbane2021component,
  title={Component-wise reduced order model lattice-type structure design},
  author={McBane, Sean and Choi, Youngsoo},
  journal={Computer methods in applied mechanics and engineering},
  volume={381},
  pages={113813},
  year={2021},
  publisher={Elsevier}
}

@article{kim2022fast,
  title={A fast and accurate physics-informed neural network reduced order model with shallow masked autoencoder},
  author={Kim, Youngkyu and Choi, Youngsoo and Widemann, David and Zohdi, Tarek},
  journal={Journal of Computational Physics},
  volume={451},
  pages={110841},
  year={2022},
  publisher={Elsevier}
}

@article{benner2015survey,
  title={A survey of projection-based model reduction methods for parametric dynamical systems},
  author={Benner, Peter and Gugercin, Serkan and Willcox, Karen},
  journal={SIAM review},
  volume={57},
  number={4},
  pages={483--531},
  year={2015},
  publisher={SIAM}
}

@article{wang2012proper,
  title={Proper orthogonal decomposition closure models for turbulent flows: a numerical comparison},
  author={Wang, Zhu and Akhtar, Imran and Borggaard, Jeff and Iliescu, Traian},
  journal={Computer Methods in Applied Mechanics and Engineering},
  volume={237},
  pages={10--26},
  year={2012},
  publisher={Elsevier}
}

@article{tsai2025local,
  title={Local Reduced-Order Modeling for Electrostatic Plasmas by Physics-Informed Solution Manifold Decomposition},
  author={Tsai, Ping-Hsuan and Chung, Seung Whan and Ghosh, Debojyoti and Loffeld, John and Choi, Youngsoo and Belof, Jonathan L},
  journal={arXiv preprint arXiv:2310.18493},
  year={2023}
}

@article{efendiev2013generalized,
  title={Generalized multiscale finite element methods (GMsFEM)},
  author={Efendiev, Yalchin and Galvis, Juan and Hou, Thomas Y},
  journal={Journal of computational physics},
  volume={251},
  pages={116--135},
  year={2013},
  publisher={Elsevier}
}

@article{calo2016randomized,
  title={Randomized oversampling for generalized multiscale finite element methods},
  author={Calo, Victor M and Efendiev, Yalchin and Galvis, Juan and Li, Guanglian},
  journal={Multiscale Modeling \& Simulation},
  volume={14},
  number={1},
  pages={482--501},
  year={2016},
  publisher={SIAM}
}

@book{maalqvist2020numerical,
  title={Numerical homogenization by localized orthogonal decomposition},
  author={Malqvist, Axel and Peterseim, Daniel},
  year={2020},
  publisher={SIAM}
}

@article{altmann2021numerical,
  title={Numerical homogenization beyond scale separation},
  author={Altmann, Robert and Henning, Patrick and Peterseim, Daniel},
  journal={Acta Numerica},
  volume={30},
  pages={1--86},
  year={2021},
  publisher={Cambridge University Press}
}

@article{chung2018constraint,
  title={Constraint energy minimizing generalized multiscale finite element method},
  author={Chung, Eric T and Efendiev, Yalchin and Leung, Wing Tat},
  journal={Computer Methods in Applied Mechanics and Engineering},
  volume={339},
  pages={298--319},
  year={2018},
  publisher={Elsevier}
}

@article{efendiev2013generalized-dg,
  title={Generalized multiscale finite element method. Symmetric interior penalty coupling},
  author={Efendiev, Yalchin and Galvis, Juan and Lazarov, R and Moon, M and Sarkis, Marcus},
  journal={Journal of Computational Physics},
  volume={255},
  pages={1--15},
  year={2013},
  publisher={Elsevier}
}

@article{cheung2020constraint,
  title={Constraint energy minimizing generalized multiscale discontinuous {Galerkin} method},
  author={Cheung, Siu Wun and Chung, Eric T and Leung, Wing Tat},
  journal={Journal of Computational and Applied Mathematics},
  volume={380},
  pages={112960},
  year={2020},
  publisher={Elsevier}
}

@article{chung2015mixed,
  title={Mixed generalized multiscale finite element methods and applications},
  author={Chung, Eric T and Efendiev, Yalchin and Lee, Chak Shing},
  journal={Multiscale Modeling \& Simulation},
  volume={13},
  number={1},
  pages={338--366},
  year={2015},
  publisher={SIAM}
}

@article{cheung2022iterative,
  title={Iterative oversampling technique for constraint energy minimizing generalized multiscale finite element method in the mixed formulation},
  author={Cheung, Siu Wun and Chung, Eric and Efendiev, Yalchin and Leung, Wing Tat and Pun, Sai-Mang},
  journal={Applied Mathematics and Computation},
  volume={415},
  pages={126622},
  year={2022},
  publisher={Elsevier}
}

@article{chung2014adaptive,
  title={An adaptive GMsFEM for high-contrast flow problems},
  author={Chung, Eric T and Efendiev, Yalchin and Li, Guanglian},
  journal={Journal of Computational Physics},
  volume={273},
  pages={54--76},
  year={2014},
  publisher={Elsevier}
}

@article{chung2017online,
  title={An online generalized multiscale discontinuous Galerkin method (GMsDGM) for flows in heterogeneous media},
  author={Chung, Eric T and Efendiev, Yalchin and Leung, Wing Tat},
  journal={Communications in Computational Physics},
  volume={21},
  number={2},
  pages={401--422},
  year={2017},
  publisher={Cambridge University Press}
}

@article{pun2023online,
  title={Online adaptive algorithm for constraint energy minimizing generalized multiscale discontinuous Galerkin method},
  author={Pun, Sai-Mang and Cheung, Siu Wun},
  journal={Multiscale Modeling \& Simulation},
  volume={21},
  number={1},
  pages={168--193},
  year={2023},
  publisher={SIAM}
}

@article{choi2025defining,
  title={Defining Foundation Models for Computational Science: A Call for Clarity and Rigor},
  author={Choi, Youngsoo and Cheung, Siu Wun and Kim, Youngkyu and Tsai, Ping-Hsuan and Diaz, Alejandro N and Zanardi, Ivan and Chung, Seung Whan and Copeland, Dylan Matthew and Kendrick, Coleman and Anderson, William and others},
  journal={arXiv preprint arXiv:2505.22904},
  year={2025}
}

@article{vasilyeva2021preconditioning,
  title={Preconditioning Markov chain Monte Carlo method for geomechanical subsidence using multiscale method and machine learning technique},
  author={Vasilyeva, Maria and Tyrylgin, Aleksei and Brown, Donald L and Mondal, Anirban},
  journal={Journal of Computational and Applied Mathematics},
  volume={392},
  pages={113420},
  year={2021},
  publisher={Elsevier}
}

\end{document}